\documentclass[reviewcopy]{elsart}
\usepackage{amsmath}%,amsthm}
\usepackage{amssymb}
\usepackage{color}
\usepackage{yfonts}
\usepackage{mathrsfs,textcomp}
\usepackage{enumerate}

\usepackage{latexsym}
\usepackage[english]{babel}
\usepackage{graphicx}
\usepackage{verbatim}
\newcommand{\SR}{{sr}}
\newcommand{\gruppi}{$H_2$, $SU(2)$, $SO(3)$, $SL(2)$ and $\mot$}
\newcommand{\mot}{SE(2)}
\newcommand{\disj}{\amalg}
\newcommand{\cont}{d\mu_\sharp}

\renewcommand{\b}[1]{{\bf #1}}
\newcommand{\Repbase}{\mathfrak{X}}
\newcommand{\dRepbase}{d\Repbase}
\newcommand{\domRbase}{{\cal H}}
\newcommand{\fzRbase}{\psi}
\newcommand{\eigRbase}{\al}

\newcommand{\ce}{\mathrm{ce}}
\newcommand{\Tanh}{\mathrm{Tanh}}
\newcommand{\Ctanh}{\mathrm{Cotanh}}
\renewcommand{\div}{\mathrm{div}}
\newcommand{\grad}{{\mathrm{grad}}}
\newcommand{\se}{\mathrm{se}}
\newcommand{\Dh}{\Delta_\SR}
\newcommand{\gradh}{\mathrm{grad}_\SR}
\newcommand{\Mat}[2]{\Pt{\ba{#1} #2\ea}}
\renewcommand{\span}[1]{\mathrm{span}\Pg{#1}}

\newcommand{\con}{{\cal C}}
\newcommand{\distr}{{\blacktriangle}}
\newcommand{\Ga}{\Gamma}
\newcommand{\Tr}[1]{\mathrm{Tr}\Pt{#1}}
\newcommand{\ad}{\mathrm{ad}\,}
\newcommand{\hp}{hypothesis}

\renewcommand{\Vec}[1]{\mathrm{Vec}\Pt{#1}}
\renewcommand{\H}{{\cal H}}
\newcommand{\g}{{\mathbf g}}
\newcommand{\metr}{\g}

\newcommand{\Hp}{{\bf (H$ _0$)}}

\newcommand{\funz}[5]{#1 : \begin{tabular}{ccl}
 #2 &$\rightarrow$& #3 \\
 #4 & $\mapsto$& #5   \end{tabular}}

\newcommand{\mmfunz}[5]{\begin{center}
\funz{$\displaystyle{#1}$}{$\displaystyle{#2}$}{$\displaystyle{#3}$}{$\displaystyle{#4}$}{$\displaystyle{#5}$}
\end{center}}
\newcommand{\fz}[3]{#1:\, #2 \rightarrow #3}

\newcommand{\bbibitem}{\bibitem}

\newcommand{\llabel}[1]{{\label{#1}}}

\newcommand{\ffoot}[1]{}
\renewcommand{\ffoot}[1]{\footnote{\noindent#1}}
\newcommand{\fffoot}[1]{}
\renewcommand{\fffoot}[1]{\footnote{\noindent#1}}

\renewcommand{\r}[1]{(\ref{#1})}

\newcommand{\bi}{\begin{itemize}}
\newcommand{\ei}{\end{itemize}}
\newcommand{\be}[1][]{\begin{enumerate}[#1]}
\newcommand{\ee}{\end{enumerate}}

\newcommand{\bd}{\begin{description}}
\newcommand{\ed}{\end{description}}
\renewcommand{\i}{\item}

\newcommand{\bqn}{\begin{eqnarray}}
\newcommand{\eqn}{\end{eqnarray}}
\newcommand{\eqnn}{\nonumber\end{eqnarray}}
\newcommand{\eqnl}[1]{\llabel{#1}\end{eqnarray}}

\newcommand{\nn}{\nonumber\\}
\newcommand{\ba}[1]{\begin{array}{#1}}
\newcommand{\ea}{\end{array}}

\newcommand{\R}{\mathbb{R}}
\newcommand{\C}{\mathbb{C}}
\newcommand{\N}{\mathbb{N}}
\newcommand{\Z}{\mathbb{Z}}

\newcommand{\proof}{{\bf Proof. }}
\newtheorem{Theorem}{\bf Theorem}
\newtheorem{lemma}[Theorem]{\bf Lemma}
\newtheorem{corollary}[Theorem]{\bf Corollary}
\newtheorem{definition}[Theorem]{\bf Definition}
\newtheorem{proposition}[Theorem]{\bf Proposition}
\newtheorem{remark}[Theorem]{\bf %\underline
{{\sl Remark}}}

\newcommand{\bt}{\begin{Theorem}}
\newcommand{\et}{\end{Theorem}}
\newcommand{\bl}{\begin{lemma}}
\newcommand{\el}{\end{lemma}}
\newcommand{\bp}{\begin{proposition}}
\newcommand{\ep}{\end{proposition}}
\newcommand{\bc}{\begin{corollary}}
\newcommand{\ec}{\end{corollary}}
\newcommand{\bdeff}{\begin{definition}}
\newcommand{\edeff}{\end{definition}}
\newcommand{\brem}{\begin{remark}\rm}
\newcommand{\erem}{\end{remark}}

\newcommand{\lam}{\lambda}
\newcommand{\al}{\alpha}

\newcommand{\ga}{\gamma}
\newcommand{\de}{\delta}

\renewcommand{\th}{\theta}

\renewcommand{\k}{{\bf k}}
\newcommand{\p}{{\mbox{\bf p}}}
\renewcommand{\l}{{\mbox{\bf L}}}
\newcommand{\Id}{\mathrm{Id}}
\newcommand{\Lie}{\mathrm{Lie}}
\newcommand{\Pt}[1]{\left( #1 \right)}
\newcommand{\Pg}[1]{\left\{ #1 \right\}}
\newcommand{\Pq}[1]{\left[ #1 \right] }
\newcommand{\Pa}[1]{\langle #1 \rangle}
\newcommand{\GFT}{GFT}
\newcommand{\DHL}{\hat{\Delta}_\SR^\lam}
\newcommand{\FF}{{\cal F}}
\renewcommand{\div}{\mbox{div}}
\newcommand{\leb}{d\mu}

\newcommand{\HS}{{\mathbf{HS}}}

\newcommand{\rrr}{}

\newcommand{\gr}[1]{{#1}}

\journal{Journal of Functional Analysis}
\begin{document}

\begin{frontmatter}
\title{The intrinsic hypoelliptic Laplacian and its  heat kernel on unimodular Lie groups}

\author{Andrei Agrachev}
\address{SISSA, via Beirut 2-4, 34014 Trieste, Italy}
\ead{agrachev@sissa.it}

\author{Ugo Boscain}
\address{LE2i,  CNRS  UMR5158, 
Universit\'e de Bourgogne,
9, avenue Alain Savary - BP 47870,   
21078 Dijon CEDEX, France}
\ead{\tt boscain@sissa.it}

\author{Jean-Paul Gauthier}
\address{Laboratoire LSIS, Universit\'e de Toulon, France}
\ead{\tt gauthier@univ-tln.fr}

\author{Francesco Rossi\corauthref{cor}}
\corauth[cor]{Corresponding author.}
\address{SISSA, via Beirut 2-4, 34014 Trieste, Italy}
\ead{\tt rossifr@sissa.it}

%%%%%%%%%%%%%%%%%%%%%%%%%%%%%%%%%%%%%%%%%%%%%%%%%%%%%%%%%%%%%%%%%%%%%%%%%%%%%%%%%
%%%%%%%%%%%%%%%%%%%%%%%%%%%%%%%%%%%%%%%%%%%%%%%%%%%%%%%%%%%%%%%%%%%%%%%%%%%%%%%%%

\newpage
\begin{abstract}
We present an invariant definition of the hypoelliptic Laplacian on sub-Riemannian structures with constant growth vector using the Popp's volume form introduced by Montgomery. This definition generalizes the one of the Laplace-Beltrami operator in Riemannian geometry. In the case of left-invariant problems on unimodular Lie groups we prove that it coincides with the usual sum of squares. 

We then extend a method (first used by Hulanicki on the Heisenberg group) to compute explicitly the kernel of the hypoelliptic heat equation on any unimodular Lie group of type I. The main tool is the noncommutative Fourier transform. We then study some relevant cases: $SU(2)$, $SO(3)$, $SL(2)$ (with the metrics inherited  by the Killing form), and the group $\mot$ of rototranslations of the plane.

%Our study is motivated by some recent results \gr{on} the cut and conjugate loci on these sub-Riemannian manifolds. \gr{Our purpose} is to understand how the singularities of the sub-Riemannian distance reflect in the kernel of the corresponding hypoelliptic heat equation.
\end{abstract}

\end{frontmatter}

\section{Introduction}

{The study of the properties of the heat kernel in a sub-Riemannian manifold  drew an increasing attention, since the \gr{pioneer} work of H\"ormander \cite{hormander}. 

Since then, many estimates and properties of the kernel in terms of the sub-Riemannian distance have been provided (see \cite{lanconelli-book,boscain-polidoro,brokko0,capogna,folland-stein,jerison,leandre,rot, varopoulos} and references therein). 

In most of the cases the hypoelliptic Laplacian appearing in the heat equation is the the sum of squares of the vector field forming an  orthonormal frame for the sub-Riemannian structure. In other cases it is built as the divergence of the horizontal gradient, where the divergence is defined using any $\con^\infty$ volume form on the manifold (see for instance \cite{taylor2}).

The Laplacians obtained in these ways are not intrinsic in the sense that they do not depend only on the sub-Riemannian distance. Indeed, 
when the Laplacian is built as the sum of squares, it depends on the choice of the orthonormal frame, while, when it is defined as divergence of the horizontal gradient, it depends on the choice of the  volume form.

The first question we address in this paper is the definition of an invariant hypoelliptic Laplacian. To our knowledge, the first time in which this question was pointed out was in a paper by 
Brockett \cite{brokko}. Many details can be found in Montgomery's book \cite{montgomery}.\\

To define the intrinsic hypoelliptic Laplacian, we proceed as in Riemannian geometry.
In Riemannian geometry  the invariant Laplacian (called the Laplace-Beltrami operator) is defined as the divergence of the gradient where the gradient is \gr{obtained} via the Riemannian metric and the divergence via the Riemannian volume form.  }

In sub-Riemannian geometry, we define the invariant hypoelliptic Laplacian as the divergence of the horizontal gradient. The horizontal gradient of a function is the natural generalization of the gradient in Riemannian geometry and it is a vector field belonging to the distribution. The divergence is computed with respect to the sub-Riemannian volume form, that can be defined for every sub-Riemannian structure with constant growth vector. This definition depends only on the sub-Riemannian structure. The sub-Riemannian volume form, called the Popp's measure, was first introduced in Montgomery's book \cite{montgomery}, where its relation with the Hausdorff measure is also discussed. The definition of the sub-Riemannian volume form is simple in the 3D contact case, and a bit more delicate in general.

We then prove that for the wide class of unimodular Lie groups (i.e. the groups where the right- and left-Haar measures coincide) the hypoelliptic Laplacian is the sum of squares for any choice of a left-invariant orthonormal base.   We recall that all compact and all nilpotent Lie groups are unimodular.\\

{In the second part of the paper, we present a method to compute explicitly the kernel of the hypoelliptic heat equation on a wide class of left-invariant sub-Riemannian structures on Lie groups. We then apply this method to the most important 3D Lie groups:
 $SU(2)$, $SO(3)$, and $SL(2)$ with the metric defined by the Killing form, the Heisenberg group $H_2$, and  the group of rototranslations of the plane $\mot$. These groups are unimodular, hence the hypoelliptic Laplacian is the sum of squares. The interest in studying
$SU(2)$, $SO(3)$ and $SL(2)$  comes from some recent results of the authors. Indeed, in \cite{nostro-gruppi} the complete description of the cut and conjugate loci for these groups was obtained. These results, together with those of this paper, open new perspectives in the direction of clarifying the relation between the presence of the cut locus and the properties of the heat kernel in the spirit of the result of Neel and Strook \cite{strook} in Riemannian geometry. Up to now the only case in which both the cut locus and the heat kernel were known explicitly was the Heisenberg group \cite{gaveau,ger,hul}.\footnote{
The Heisenberg group is in a sense a very  degenerate example. For instance, in this case the cut locus coincides globally  with the conjugate locus (set of points where geodesics lose local optimality) and  many properties that one expects to be  distinct for more generic situations cannot be distinguished. The application of our method 
to the Heisenberg group $H_2$ provides in a few lines the Gaveau-Hulanicki formula \cite{gaveau,hul}.
} 
%(The hypoelliptic heat kernel is however known for all step 2 nilpotent Lie groups \cite{greiner}.)

The interest in the hypoelliptic heat kernel on $\mot$ comes from  a model of human vision. It was recognized in 
\cite{citti-sarti,petitot} that  the visual cortex V1 solves a nonisotropic diffusion problem on the group $SE(2)$ while reconstructing a partially hidden or corrupted image. The study of the cut locus on $SE(2)$ is a work in progress. Preliminary results can be found in \cite{igorino-yuri}. }

%{\color{red} The method we present is general and it is based on the generalized Fourier transform (\GFT, for short): it is the transformation of a function from a Lie group $G$ to $\C$ into a function from the dual space $\hat{G}$ (i.e. the set of all non-equivalent irreducible representations) to Hilbert-Schmidt operators. } 
The method is based upon the generalized (noncommutative) Fourier transform (\GFT, for short), that permits to disintegrate\footnote{One could also say decompose (possibly continuously).} a function from a Lie  group $G$ to $\R$ on its components on (the class of) non-equivalent unitary irreducible representations of $G$. 
This technique permits to transform the hypoelliptic heat equation into an equation in  the dual of the group\footnote{\label{nota-dual} In this paper, by the dual of the group, we mean  the support of the Plancherel measure on the set of non-equivalent unitary irreducible representations of $G$; we thus ignore the singular representations.}, that is particularly simple since the GFT disintegrate the right-regular representations
and the hypoelliptic Laplacian is built with left-invariant vector fields (to which a one parameter group of right-translations is associated).

Unless we are in the abelian case, the dual of a Lie group in general is not a group. In the compact case it is a so called Tannaka category \cite{hewittII,tannaka}
and it is a discrete set. In the nilpotent case it has the structure of $\R^n$ for some $n$. In the general case it can have a quite complicated structure.  However, under certain hypotheses (see Section \ref{s-genfour}), it is a measure space if endowed with the so called Plancherel measure. Roughly speaking, the \GFT\ is an isometry between $L^2(G,\C)$ (the set of complex-valued square integrable functions over $G$, with respect to the Haar measure) and the set of Hilbert-Schmidt operators with respect to the Plancherel measure.

The difficulties of applying our method in specific cases rely mostly on two points:
\be[i)]
\i computing the tools for the \GFT, i.e. the non-equivalent irreducible representations of the group and the Plancherel measure. This is a difficult problem in general: however, for certain classes of Lie groups there are suitable techniques  (for instance the Kirillov orbit method for nilpotent Lie groups \cite{kirillov}, or methods for semidirect products).  For the groups \gr{disccussed} in this paper, the sets of non-equivalent irreducible representations (and hence the \GFT) are well known (see for instance \cite{sugiura});
\i finding the spectrum of an operator (the \GFT\ of the hypoelliptic Laplacian). Depending on the structure of the group and on its dimension, this problem gives rise to a matrix equation, an ODE or a PDE.
\ee
Then one can express the kernel of the hypoelliptic heat equation in terms of eigenfunctions of the \GFT\ of the hypoelliptic Laplacian, or in terms of the kernel of the transformed equation.

For the cases treated in this paper, we have the following (the symbol $\disj$ means disjoint union):
\begin{center}
\begin{tabular}{|l|l|l|l|}\hline
Group& Dual of  & \GFT\ &Eigenfunctions of the \GFT\\
&the group& of the hypoelliptic Laplacian& of the hypoelliptic Laplacian\\
\hline
$H_2$& $\R$ & $\displaystyle \frac{d^2}{dx^2} -\lam^2 x^2     $ (quantum Harmonic oscillator)&Hermite polynomials\\\hline
% \vspace{-2mm} &&&\\\hline
 $SU(2)$& $\N$ & Linear finite dimensional operator  related to the& Complex homogeneous  \\
 && quantum angular momentum& polynomials in two variables\\\hline
 $SO(3)$& $\N$ &  Linear finite dimensional operator   related to& Spherical harmonics\\ 
 && orbital quantum angular momentum&\\\hline
 $SL(2)$ & $\R^+\disj\R^+$ & Continuous: Linear operator on analytic functions &Complex monomials\\
&$\disj\N\disj\N$&~~~~~~~~~with domain $\Pg{|x|=1}\subset\C$  &\\
&&Discrete: Linear operator on analytic functions&\\
&&~~~~~~~~~with domain $\Pg{|x|<1}\subset\C$&
\\\hline
 $\mot$& $\R^+$  &  $\displaystyle \frac{d^2}{d\th^2} -\lam^2 \cos^2(\th)     $ (Mathieu's equation)& $2\pi$-periodic Mathieu functions  \\\hline\hline 

\end{tabular}
\end{center}

The idea of using the \GFT\ to compute the hypoelliptic heat kernel is not new: it was already used on the Heisenberg group in \cite{hul} at the same time as the Gaveau formula was published in \cite{gaveau}, and on all step 2 nilpotent Lie groups in \cite{cyl,greiner}. See also the related work \cite{klingler}.

The structure of the paper is the following: in Section \ref{s-Dh} we recall some basic definitions from sub-Riemannian geometry and we construct the sub-Riemannian volume form. We then give the definition of the hypoelliptic Laplacian on a regular sub-Riemannian manifold, and we show that the hypothesis of regularity cannot be dropped in general. For this purpose, we show that the invariant hypoelliptic Laplacian defined on the Martinet sub-Riemannian structure is singular. We then move to left-invariant sub-Riemannian structures on Lie groups and we show that a Lie group is unimodular if and only if the invariant hypoelliptic Laplacian is the sum of squares.
We also provide an example of a 3D non-unimodular Lie group for which the invariant hypoelliptic Laplacian is not the sum of squares.
The section ends with the proof that the invariant hypoelliptic Laplacian can be expressed as
$$  \Dh=-\sum_{i=1}^m L_{X_i}^\ast   L_{X_i}, $$
where the formal adjoint $L_{X_i}^\ast$ is built with the sub-Riemannian volume form, providing a connection with existing literature (see e.g. \cite{jerison2}). The invariant hypoelliptic Laplacian is then the sum of squares when $L_{X_i}$ are skew-adjoint\footnote{This point of view permits to give an alternative proof of the fact that the invariant hypoelliptic Laplacian on left-invariant structures on unimodular Lie groups is the sum of squares. \gr{As a matter of fact}, left-invariant vector fields are formally skew-adjoint with respect to the right-Haar measure.
On Lie groups the invariant volume form is left-invariant, hence is proportional to the left-Haar measure, and is in turn proportional to the right-Haar measure on unimodular groups.}.

In Section \ref{s-genfour} we recall basic tools of the \GFT\ and we describe our general method to compute the heat kernel of the hypoelliptic Laplacian on unimodular Lie groups of type I. We provide two useful formulas, one in the case where the GFT of the hypoelliptic Laplacian has discrete spectrum, and the other in the case where the GFT of the hypoelliptic heat equation admits a kernel. 

In Section \ref{s-examples} we apply our method to compute the kernel on \gruppi. For the Heisenberg group we use the formula involving the kernel of the transformed equation (the Mehler kernel). For the other groups we use the formula in terms of eigenvalues and eigenvectors of the GFT of the hypoelliptic Laplacian.

The application of our method to higher dimensional sub-Riemannian problems and in particular to the nilpotent Lie groups $(2,3,4)$ (the Engel group) and $(2,3,5)$ is the subject of a forthcoming paper.

\section{The hypoelliptic Laplacian}
\llabel{s-Dh}
In this Section we give a definition of the hypoelliptic Laplacian $\Dh$ on a regular sub-Riemannian manifold $M$.

\subsection{Sub-Riemannian manifolds}
We start recalling the definition of sub-Riemannian manifold.
\bdeff
A $(n,m)$-sub-Riemannian manifold is a triple $(M,\distr,{\mathbf g})$, 
where
\bi
\i $M$ is a connected smooth manifold of dimension $n$;
\i $\distr$ is a smooth distribution of constant rank $m\leq n$ satisfying the {\bf H\"ormander condition}, i.e. $\distr$ is a smooth map that associates to $q\in M$  a $m$-dim subspace $\distr(q)$ of $T_qM$ and $\forall~q\in M$ we have
\bqn\llabel{Hor}\span{[X_1,[\ldots[X_{k-1},X_k]\ldots]](q)~|~X_i\in\mathrm{Vec}_H(M)}=T_qM
\eqn
where $\mathrm{Vec}_H(M)$ denotes the set of {\bf horizontal smooth vector fields} on $M$, i.e. $$\mathrm{Vec}_H(M)=\Pg{X\in\mathrm{Vec}(M)\ |\ X(p)\in\distr(p)~\ \forall~p\in M}.$$
\i $\g_q$ is a Riemannian metric on $\distr(q)$, that is smooth 
as function of $q$.
\ei
When $M$ is an orientable manifold, we say that the sub-Riemannian manifold is orientable.
\edeff

\brem
Usually sub-Riemannian manifolds are defined with $m<n$. In our definition we decided to include the Riemannian case $m=n$, since all our results hold  in that case. Notice that if $m=n$ then condition \r{Hor} is automatically satisfied.
\erem

A Lipschitz continuous curve $\ga:[0,T]\to M$ is said to be \b{horizontal} if 
$\dot\ga(t)\in\distr(\ga(t))$ for almost every $t\in[0,T]$.

Given an horizontal curve $\ga:[0,T]\to M$, the {\it length of $\ga$} is
\bqn
l(\ga)=\int_0^T \sqrt{ \g_{\ga(t)} (\dot \ga(t),\dot \ga(t))}~dt.
\eqnl{e-lunghezza}
The {\it distance} induced by the sub-Riemannian structure on $M$ is the 
function
\bqn
d(q_0,q_1)=\inf \{l(\ga)\mid \ga(0)=q_0,\ga(T)=q_1, \ga\ \mathrm{horizontal}\}.
\eqnl{e-dipoi}

The \hp\ of connectedness of M and the H\"ormander condition guarantee the finiteness and the continuity of $d(\cdot,\cdot)$ with respect to the topology of $M$ (Chow's Theorem, see for instance \cite{agra-book}). The function $d(\cdot,\cdot)$ is called the Carnot-Charateodory distance and gives to $M$ the structure of metric space (see \cite{bellaiche,gromov}).

It is a standard fact that $l(\ga)$ is invariant under reparameterization of the curve $\ga$.
Moreover, if an admissible curve $\ga$ minimizes the so-called {\it energy functional}
$$ E(\ga)=\int_0^T {\g}_{\ga(t)}(\dot \ga(t),\dot \ga(t))~dt. $$
with $T$ fixed (and fixed initial and final point), then $v=\sqrt{\g_{\ga(t)}(\dot \ga(t),\dot \ga(t))}$
is constant and $\ga$ is also a minimizer of $l(\cdot)$.
On the other side, a minimizer $\ga$ of $l(\cdot)$ such that  $v$ is constant is a minimizer of $E(\cdot)$ with $T=l(\ga)/v$.

A {\it geodesic} for  the sub-Riemannian manifold  is a curve $\ga:[0,T]\to M$ such that for every sufficiently small interval $[t_1,t_2]\subset [0,T]$, $\ga_{|_{[t_1,t_2]}}$ is a minimizer of $E(\cdot)$.
A geodesic for which $\g_{\ga(t)}(\dot \ga(t),\dot \ga(t))$  is (constantly) equal to one is said to be parameterized by arclength.

Locally, the pair $(\distr,{\mathbf g})$ can be given by assigning a set of $m$ smooth vector fields spanning $\distr$ and that are orthonormal for ${\mathbf g}$, i.e.  
\bqn
\distr(q)=\span{X_1(q),\dots,X_m(q)}, ~~~\metr_q(X_i(q),X_j(q))=\delta_{ij}.
\eqnl{trivializable}
In this case, the set $\Pg{X_1,\ldots,X_m}$ is called a local {\bf orthonormal frame} for the sub-Riemannian structure. When  $(\distr,{\mathbf g})$ can be defined as in \r{trivializable} by $m$ vector fields defined globally, we say that the sub-Riemannian manifold is {\it trivializable}. 

%%%%%%%%%%%%%%%%%%%%%%%%%%%%%%%%%%%%%%%%%%%%%%%%%%%%%%%%%%%%%%%%%%%%

Given a $(n,m)$- trivializable sub-Riemannian manifold, the problem of finding a curve minimizing the energy between two fixed points  $q_0,q_1\in M$ is
naturally formulated as the optimal control problem 
\bqn
\dot q(t)&=&\sum_{i=1}^m u_i(t) X_i(q(t))\,,~~~u_i(.)\in L^\infty([0,T],\R)\,,
~~~\int_0^T \sum_{i=1}^m u_i^2(t)~dt\to\min,\label{sopra}\\ q(0)&=&q_0,~~~q(T)=q_1.
\eqnn

It is a standard fact that this optimal control problem is equivalent to the minimum time problem with controls $u_1,\ldots, u_m$ satisfying $u_1(t)^2+\ldots+u_m(t)^2\leq 1$ in $[0,T]$.

When the manifold is analytic and the orthonormal frame can be assigned through $m$ analytic vector fields, we say that the sub-Riemannian manifold is {\it analytic}.

We end this section with the definition of the small flag of the distribution $\distr$:
\bdeff
Let $\distr$ be a distribution and define through the recursive formula
$$\distr_1:=\distr~~~\distr_{n+1}:=\distr_n+[\distr_n,\distr]$$
where $\distr_{n+1}(q_0):=\distr_n(q_0)+[\distr_n(q_0),\distr(q_0)]=\\
=\Pg{X_1(q_0)+[X_2,X_3](q_0)\ |\ X_1(q),X_2(q)\in\distr_n(q),\ X_3(q)\in\distr(q)~~\forall~q\in M}.$ The small flag of $\distr$ is the sequence 
$$\distr_1\subset\distr_2\subset\ldots\subset\distr_n\subset\ldots$$

A sub-Riemannian manifold is said to be {\bf regular} if for each $n=1,2,\ldots$ the dimension of $\distr_n(q_0)=\Pg{f(q_0)\ |\ f(q)\in\distr_n(q)~\forall\ q\in M}$ does not depend on the point $q_0\in M$.

A 3D sub-Riemannian manifold is said to be a {\bf 3D contact manifold} if $\distr$ has dimension 2 and $\distr_2(q_0)=T_{q_0}M$ for any point $q_0\in M$.
\edeff
In this paper we always deal with regular sub-Riemannian manifolds.

\subsubsection{Left-invariant sub-Riemannian manifolds}
\label{ss-leftmanifold}

In this section we present a natural sub-Riemannian structure that can be defined on Lie groups. All along the paper, we use the notation for Lie groups of matrices. For general Lie groups, by $gv$ with $g\in G$ and $v\in \l$, we mean $(L_g)_*(v)$ where $L_g$ is the left-translation of the group.

\bdeff 
Let $G$ be a Lie group with Lie algebra $\l$ and $\p\subseteq\l$ a subspace of $\l$ satisfying the {\bf Lie bracket generating condition} $$\Lie~\p:=\span{[p_1,[p_2,\ldots,[p_{n-1},p_n]]]\ |\ p_i\in\p}=\l.$$
Endow $\p$ with a positive definite quadratic form $\Pa{.,.}$. Define a sub-Riemannian structure on $G$ as follows:
\bi
\i the distribution is the left-invariant distribution $\distr(g):=g\p$;
\i the quadratic form $\metr$ on $\distr$ is given by $\metr_g(v_1,v_2):=\Pa{g^{-1}v_1,g^{-1}v_2}$.
\ei
In this case we say that $(G, \distr, \metr)$ is a left-invariant sub-Riemannian manifold.
\llabel{d-lieg-leftinv}
\edeff
\brem Observe that all left-invariant manifolds $(G, \distr, \metr)$ are regular.
\erem

In the following we define a left-invariant sub-Riemannian manifold choosing a set of $m$ vectors $\Pg{p_1,\ldots,p_m}$ being an orthonormal basis for the subspace $\p\subseteq\l$ with respect to the metric defined in Definition \ref{d-lieg-leftinv}, i.e. $\p=\span{p_1,\ldots,p_m}$ and $\Pa{p_i,p_j}=\de_{ij}$. We thus have $\distr(g)=g\p=\span{gp_1,\ldots,gp_m}$ and $\g_g(gp_i,gp_j)=\de_{ij}$. Hence every left-invariant sub-Riemannian manifold is trivializable.

The problem of finding the minimal energy between the identity and a point $g_1\in G$ in fixed time $T$ becomes the left-invariant optimal control problem
\bqn
\dot g(t)&=&g(t)\left(\sum_i u_i(t) p_i\right),~~~~~~~~u_i(.)\in L^\infty([0,T],\R)
~~~~~~~\int_{0}^{T}\sum_i u_i^2(t)~dt\to\min,\label{e-controllo}\\
g(0)&=&\Id,~~~g(T)=g_1.
\eqnn
\brem \llabel{rem-minesiste}
This problem admits a solution, see for instance Chapter 5 of \cite{piccoli}.
\erem

\subsection{Definition of the hypoelliptic Laplacian on a sub-Riemannian manifold}
\newcommand{\muh}{\mu_\SR}
\newcommand{\divh}{\mathrm{div}_\SR}

In this section we define the intrinsic hypoelliptic Laplacian on a regular orientable sub-Riemannian manifold $(M,\distr,\metr)$. 
This definition generalizes the one of the  Laplace-Beltrami operator on an orientable Riemannian manifold, that is $\Delta \phi:=\div~\grad~\phi$, where $\grad$ is the unique operator from  $\mathcal{C}^\infty(M)$ to $\Vec{M}$ satisfying $\g_q(\grad~\phi(q),v)=d\phi_q (v)~~~~\forall~q\in M,~v\in T_qM,$ and the divergence of a vector field $X$ is the unique function satisfying   
$\div X \mu=L_X \mu$ where $\mu$ is the Riemannian volume form.

We first define the sub-Riemannian gradient of a function, that is an horizontal vector field.
\bdeff
\llabel{def-gradh}
Let $(M,\distr,\g)$ be a sub-Riemannian manifold: the {\bf horizontal gradient} is the unique operator $\gradh$ from $\mathcal{C}^\infty(M)$ to ${{\mathrm{Vec}}_H(M)}$ 
satisfying $\g_q(\gradh\phi(q),v)=d\phi_q (v)~~~~\forall~q\in M,~v\in \distr(q).$
\edeff
One can easily check that if $\{X_1,\ldots X_m \}$ is a local orthonormal frame for $(M,\distr,\metr)$, then $\gradh \phi=\sum_{i=1}^m \Pt{L_{X_i}\phi}X_i$.

The question of defining a sub-Riemannian volume form is more delicate. We start by considering the 3D contact case. 
\bp
Let  $(M,\distr,\g)$ be an orientable  3D contact  sub-Riemannian structure and $\{X_1,X_2 \}$  a local orthonormal frame. Let $X_3=[X_1,X_2]$ and $dX_1$, $dX_2$, $dX_3$ the dual basis, i.e. $dX_i(X_j)=\de_{ij}$. Then $\muh:=dX_1\wedge dX_2 \wedge dX_3$ is an intrinsic volume form, i.e. it is invariant for a orientation preserving change of orthonormal frame.
\ep
\proof Consider two different orthonormal frames with the same orientation $\Pg{X_1,X_2}$ and $\Pg{Y_1,Y_2}$. We have to prove that $dX_1\wedge dX_2
\wedge dX_3=dY_1\wedge dY_2 \wedge dY_3$ with $X_3=[X_1,X_2]$,
$Y_3=[Y_1,Y_2]$. We have
$$\Mat{c}{Y_1\\Y_2}=\Mat{cc}{\cos(f(q))&\sin(f(q))\\-\sin(f(q))&\cos(f(q))}
\Mat{c}{X_1\\X_2},$$
for some real-valued smooth function $f$. A direct computation shows that 
\bqn
Y_3=X_3+f_1X_1+f_2X_2
\eqnl{e-X3Y3} where $f_1$ and $f_2$ are two smooth functions depending on $f$.

We first prove that $dX_1\wedge dX_2=dY_1\wedge dY_2$. Since 
the change of
variables $\Pg{X_1,X_2}\mapsto\Pg{Y_1,Y_2}$ is norm-preserving, 
we have $dX_1\wedge
dX_2(v,w)=dY_1\wedge dY_2(v,w)$ when $v,w\in\distr$. 
Consider now any 
vector $v=v_1X_1+v_2X_2+v_3X_3=v_1'Y_1+v_2'Y_2+v_3'Y_3$: as a consequence of
\r{e-X3Y3}, we have $v_3=v_3'$. Take another vector
$w=w_1X_1+w_2X_2+w_3X_3=w_1'Y_1+w_2'Y_2+w_3Y_3$ and compute $$dX_1\wedge
dX_2(v,w)=dX_1\wedge dX_2(v-v_3X_3,w-w_3X_3)=dY_1\wedge
dY_2(v-v_3X_3,w-w_3X_3)=dY_1\wedge dY_2(v,w),$$ because the vectors
$v-v_3X_3,w-w_3X_3$ are horizontal. Hence the two 2-forms coincide.

From \r{e-X3Y3} we also have $dY_3=dX_3+f_1'dX_1+f_2'dX_2$ for some smooth functions $f'_1,f'_2$. Hence 
we have $dY_1\wedge dY_2\wedge dY_3=dX_1\wedge
dX_2\wedge dY_3=dX_1\wedge dX_2\wedge (dX_3+f_1'dX_1+f_2'dX_2)=dX_1\wedge
dX_2\wedge dX_3,$ where the last identity is a consequence of skew-symmetry
of differential forms. 
\qed

\brem
Indeed,  even if in the 3D contact case there is no  scalar product in $T_qM$, it is possible to define a  natural volume form, since on $\distr$ the scalar product is defined by $\g$ and formula \r{e-X3Y3} guarantees the existence of a natural scalar product in  $(\distr+[\distr,\distr])/\distr$.
\erem

The previous result generalizes to any regular orientable sub-Riemannian structure, as presented below.

\subsubsection{Definition of the intrinsic volume form}

Let $0=E_0\subset E_1\subset\ldots\subset E_k=E$ be a filtration
of an $n$-dimensional vector space $E$. Let $e_1,\ldots,e_n$ be a basis
of $E$ such that $E_i=\span{e_1,\ldots,e_{n_i}}$. Obviously, the
wedge product $ e_1\wedge\ldots\wedge e_n $ depends only on the
residue classes
$$\bar e_j=(e_j+E_{i_j})\in E_{i_j+1}/E_{i_j},$$
where $n_{i_j}<j\le n_{i_j+1},\ j=1,\ldots,n$. This property
induces a natural (i\,e. independent on the choice of the basis)
isomorphism of 1-dimensional spaces:
$$\bigwedge\nolimits^nE\cong\bigwedge^n\left(\bigoplus\limits_{i=1}^k(E_i/E_{i-1})\right).
$$

Now consider the filtration
$$
0\subset\distr_1(q)\subset\ldots\subset\distr_k(q)=T_qM,\quad\dim\distr_i(q)=n_i.
$$
Let $X_1,\ldots,X_i$ be smooth sections of $\distr=\distr_1$; then
the vector
$$
\bigl([X_1,[\ldots,X_i]\ldots](q)+\distr_{i-1}(q)\bigr)\in\distr_i(q)/\distr_{i-1}(q)
$$
depends only on $X_1(q)\otimes\ldots\otimes X_i(q)$.

We thus obtain a well-defined surjective linear map
\mmfunz{\beta_i}{\distr(q)^{\otimes i}}{\distr_i(q)/\distr_{i-1}(q)}{X_1(q)\otimes\ldots\otimes X_i(q)}{\bigl([X_1,[\ldots,X_i]\ldots](q)+\distr_{i-1}(q)\bigr)}

The Euclidean structure on $\distr(q)$ induces an Euclidean structure
on $\distr(q)^{\otimes i}$ by the standard formula:
$$
\langle\xi_1\otimes\ldots\otimes\xi_i,\eta_1\otimes\ldots\otimes\eta_i\rangle=
\langle\xi_1,\eta_1\rangle\ldots\langle\xi_i,\eta_i\rangle,\quad
\xi_j,\eta_j\in\distr(q),\ j=1,\ldots,i.
$$
Then the formula: $$ |v|=\min\{|\bar\xi|:
\bar\xi\in\beta^{-1}_i(v)\},\quad v\in\distr^i(q)/\distr^{i-1}(q)
$$
defines an Euclidean norm on $\distr^i(q)/\distr^{i-1}(q)$.

Let $\nu_i$ be the volume form on $\distr^i(q)/\distr^{i-1}(q)$
associated with the Euclidean structure:
$$
\langle\nu_i,v_1\wedge\ldots v_{m_i}\rangle= \det\nolimits^{\frac
1{m_i}}\left\{\langle v_j,v_{j'}\rangle\right\}_{j,j'=1}^{m_i},
$$
where $m_i=n_i-n_{i-1}=\dim(\distr^i(q)/\distr^{i-1}(q))$.

Finally, the intrinsic volume form $\muh$ on $T_qM$ is the image of
$\nu_1\wedge\ldots\wedge\nu_k$ under the natural isomorphism
$$
\bigwedge\nolimits^n\left(\bigoplus\limits_{i=1}^k(\distr^i(q)/\distr^{i-1}(q))\right)^*
\cong\bigwedge\nolimits^n\left(T_qM\right)^*.
$$

\brem The construction given above appeared for the first time in the book of Mongomery \cite[Section 10.5]{montgomery}. He called the measure $\muh$ the Popp's measure. He  also observed that a sub-Riemannian volume form was the  only missing ingredient to get an intrinsic definition of hypoelliptic Laplacian.\footnote{{\rrr 
Montgomery did not use Popp's measure to get the intrinsic definition of the hypoelliptic Laplacian since there 
are, a priori, two natural measures on a regular sub-Riemannian manifold: the Popp's measure, and the 
Hausdorff measure (see \cite{mitchell, pansu}). However, for left-invariant sub-Riemannian manifolds,  both of them are proportional to the left-Haar measure. See Remark \ref{r-misure}.} }
%However Montgomery did not applied the Popp's measure to find an explicit expression of the invariant hypoelliptic Laplacian, as for instance formula \r{f-esplicita} below.
\erem
%\brem
%\label{r-hausdorff} 
%One could argue that an intrinsic volume form in sub-Riemannian geometry was already available, namely the Hausdorff measure. However the Popp's measure has two advantages with respect to the Hausdorff measure: {\bf i)} it can be easily computed,  {\bf ii)} it has smooth density. 

%On the other side, the Hausdorff measure is very difficult to be computed in practice and, even if it is absolutely continuous with respect to the Lebesgue measure (\cite{mitchell,pansu}), it is not known if it is smooth. A natural question is therefore (see \cite[p. 141]{montgomery}):\\\\
%
%{\bf Q2}: are the Popp's measure and the Hausdorff measure proportional?\\\\
%
%For left-invariant structures on Lie groups, both  measures are left-invariant, hence they are proportional to  the left-Haar measure (see below for the definition), however in the general case the answer is unknown.
%\erem
Once the volume form is defined, the divergence of a vector field $X$ is defined as in Riemannian geometry, i.e. it is the function $\divh X$ satisfying
$\divh X \muh=L_X\muh$. We are now ready to define the intrinsic hypoelliptic Laplacian.
\bdeff
Let  $(M,\distr,\g)$ be an orientable  regular sub-Riemannian manifold. Then the intrinsic hypoelliptic Laplacian is
$\Dh\phi:=\divh\gradh\phi$.
\edeff

Consider now an orientable regular sub-Riemannian structure $(M,\distr,\g)$ and let $\{X_1,\ldots X_m \}$ be a local orthonormal frame. We want to find an explicit expression for the operator $\Dh$.  If $n=m$ then $\Dh$ is the Laplace Beltrami operator. Otherwise 
consider $n-m$ vector fields $X_{m+1},\ldots,X_n$ such that $\Pg{X_1(q),\ldots,X_m(q),X_{m+1}(q),\ldots,X_n(q)}$ is a basis of $T_qM$ for all $q$ in a certain open set $U$. The volume form $\muh$ is $\muh=f(q)dX_1\wedge\ldots\wedge dX_n$, with $dX_i$ dual basis of $X_1,\ldots,X_n$: then we can find other $n-m$ vector fields, that we still call  $X_{m+1},\ldots,X_n$, for which we have $\muh=dX_1\wedge\ldots\wedge dX_n$.

Recall that $\Dh\phi$ satisfies $(\Dh\phi)\muh=L_X\muh$ with $X=\gradh \phi$. We have 
\bqn
L_X \muh&=&\sum_{i=1}^m (-1)^{i+1} \left\{ d\Pt{\Pa{d\phi,X_i}}\wedge dX_1\wedge\ldots\wedge\hat{dX_i}\wedge\ldots\wedge dX_n +\right.\nn
&&\left.+\Pa{d\phi,X_i} d\Pt{dX_1\wedge\ldots\wedge\hat{dX_i}\wedge\ldots\wedge dX_n}\right\}.
\eqnn
Applying standard results of differential calculus, we have $d\Pt{\Pa{d\phi,X_i}}\wedge dX_1\wedge\ldots\wedge\hat{dX_i}\wedge\ldots\wedge dX_n=(-1)^{i+1} L_{X_i}^2\phi~ \muh$ and $d\Pt{dX_1\wedge\ldots\wedge\hat{dX_i}\wedge\ldots\wedge dX_n}=(-1)^{i+1} \Tr{\ad X_i}\muh$, where the adjoint map is  
\mmfunz{\ad X_i}{\Vec{U}}{\Vec{U}}{X}{[X_i,X].} and by $\Tr{\ad X_i}$ we mean $\sum_{j=1}^n dX_j([X_i,X_j])$. Finally, we find the expression 
\bqn
\label{f-esplicita}
\Dh\phi=\sum_{i=1}^m \Pt{L_{X_i}^2\phi + L_{X_i}\phi~ \Tr{\ad X_i}}.
\eqn
Notice that the formula depends on the choice of the vector fields $X_{m+1},\ldots,X_n$.

The hypoellipticity of $\Dh$ (i.e. given $U\subset M$ and $\fz{\phi}{U}{\C}$ such that $\Dh\phi\in\con^\infty$, then $\phi$ is $\con^\infty$) follows from the H\"ormander Theorem (see \cite{hormander}):
\bt
Let $L$ be a differential operator on a manifold $M$, that locally in a neighborhood $U$ is written as 
$L=\sum_{i=1}^m L_{X_i}^2+L_{X_0}$, where $X_0, X_1\ldots,X_m$ are $\con^\infty$ vector fields. If   
$\Lie_q\{X_0,X_1,\ldots,X_m \}=T_q M$ for all $q\in U$, then $L$ is hypoelliptic. 
\et
Indeed, $\Dh$ is written locally as $\sum_{i=1}^m L_{X_i}^2+L_{X_0}$ with the first-order term $L_{X_0}=\sum_{i=1}^m\Tr{\ad X_i}~L_{X_i}$. Moreover by hypothesis we have that $\Lie_q\{X_1,\ldots,X_m \}=T_q M$, hence the  H\"ormander theorem applies.\\

\brem
Notice that in the Riemannian case, i.e. for $m=n$, $\Dh$ coincides with the Laplace-Beltrami operator.
\erem
\brem
The hypothesis that the sub-Riemannian manifold is regular is crucial for the construction of the invariant volume form. For instance for the Martinet metric on $\R^3$, that is the sub-Riemannian structure  for which  
$L_1=\partial_x+\frac{y^2}{2}\partial_z$ and 
$L_2=\partial_y$ form an orthonormal base, one gets on $\R^3\setminus\Pg{y=0}$
$$
%\muh=\frac1ydx\wedge dy\wedge dz \mbox{ ~~and~~ }
\Dh= (L_1)^2+(L_2)^2-\frac1y L_2.
$$
This is not surprising at all. \gr{As a matter of fact}, even the Laplace-Beltrami operator is singular  in almost-Riemannian geometry (see \cite{gauss-bonnet} and references therein). For instance, for the Grushin metric on $\R^2$, that is the singular Riemannian structure   for which $L_1=\partial_x$ and $L_2=x \,\partial_y$ form an orthonormal frame, one gets on $\R^2\setminus\{x=0 \}$
$$
\Delta_{L.B.}= (L_1)^2+(L_2)^2-\frac1x L_1.
$$
\erem

\subsection{The hypoelliptic Laplacian on Lie groups}
In the case of left-invariant sub-Riemannian manifolds, there is an intrinsic  global expression of $\Dh$.
\bc
\label{cor-c}
Let $(G,\distr,\metr)$ be a left-invariant sub-Riemannian manifold generated by the orthonormal basis $\Pg{p_1,\ldots,p_m}\subset\l$. Then the hypoelliptic Laplacian is
\bqn
\Dh\phi&=&\sum_{i=1}^m \Pt{L_{X_i}^2\phi + L_{X_i}\phi~ \Tr{\ad p_i}}\llabel{e-Dh-Lie}
\eqn
where $L_{X_i}$ is the Lie derivative w.r.t. the field $X_i=gp_i$.
\ec
\proof
If $m\leq n$, we can find $n-m$ vectors $\Pg{p_{m+1},\ldots,p_n}$ such that $\Pg{p_1,\ldots,p_n}$ is a basis for $\l$. Choose the fields $X_i:=gp_i$ and follow the computation given above: we find formula \r{e-Dh-Lie}. In this case the adjoint map is intrinsically defined and the trace does not depend on the choice of $X_{m+1},\ldots,X_{n}$.
\qed\\\\
\newcommand{\fmod}{\Psi}
\newcommand{\Ad}{\mathrm{Ad}}
The formula above reduces to the sum of squares in the wide class of unimodular Lie groups. We recall that on a Lie group of dimension $n$, there always exist  a left-invariant $n$-form $\mu_L$ and a right-invariant  $n$-form $\mu_R$ (called respectively left- and right-Haar measures), that are unique up to a multiplicative constant.
These forms have the properties that 
$$
\int_G f(a g)\mu_L(g)=\int_G f(g)\mu_L(g),~~~~~\int_G f(g a)\mu_R(g)=\int_G f(g)\mu_R(g),\mbox{ for every  }f\in L^1(G,\R)\mbox{ and }a\in G, 
$$
where $L^1$ is intended with respect to the left-Haar measure in the first identity and with respect to the right-Haar measure in the second one. The group is called unimodular if $\mu_L$ and $\mu_R$ are proportional. 
\brem
\label{r-misure}
Notice that for left-invariant sub-Riemannian manifolds the intrinsic volume form and the Hausdorff measure $\mu_H$ are left-invariant, hence they are proportional to the the left Haar measure $\mu_L$. On unimodular Lie groups one can assume $\muh=\mu_L=\mu_R=\al\mu_H$, where $\al>0$ is a constant that is unknown even for the simplest \gr{among the} genuine sub-Riemannian structures, i.e. the Heisenberg group.
\erem

\bp
Under the hypotheses of Corollary \ref{cor-c}, if $G$ is unimodular then 
$\Dh\phi=\sum_{i=1}^m L_{X_i}^2\phi.$
\ep
\proof
Consider the modular function $\fmod$, that is the unique function such that $\int_G f(h^{-1}g)\mu_R(g)=\fmod(h)\int_G f(g)\mu_R(g)$ for all $f$ measurable. It is well known that $\fmod(g)=\det(\Ad_g)$ and that $\fmod(g)\equiv 1$ if and only if $G$ is unimodular.

Consider a curve $\ga(t)$ such that $\dot{\ga}$ exists for $t=t_0$: then $\ga(t)=g_0e^{(t-t_0)\eta+o(t-t_0)}$ with $g_0=\ga(t_0)$ and for some $\eta\in\l$. We have
\bqn
\frac{d}{dt}_{|_{t=t_0}}\det(\Ad_{\ga(t)})&=&\Tr{(\Ad_{g_0})^{-1} \Pq{\frac{d}{ds}_{|_{s=0}} \Ad_{g_0 e^{s\eta+o(s)}}}}=\label{e-derAd}\\
&=&\Tr{\Ad_{g_0^{-1}}\Ad_{g_0} \ad_\eta}=\Tr{\ad_\eta}.
\eqnn
Now choose the curve $\ga(t)=g_0e^{tp_i}$ and observe that $\det(\Ad_{\ga(t)})\equiv 1$, then $\Tr{\ad_{p_i}}=0$. The conclusion follows from \r{e-Dh-Lie}.
\qed\\\\
All  the groups treated in this paper, (i.e. \gruppi) are unimodular.  Hence the invariant hypoelliptic Laplacian is the sum of squares. A kind of inverse result holds:
\bp
\label{p-inversa}
Let $(G,\distr,\metr)$ be a left-invariant sub-Riemannian manifold generated by the orthonormal basis $\Pg{p_1,\ldots,p_m}\subset\l$. If the hypoelliptic Laplacian satisfies $\Dh\phi=\sum_{i=1}^m L_{X_i}^2\phi$, then $G$ is unimodular.
\ep
\proof
We start observing that $\Dh\phi=\sum_{i=1}^m L_{X_i}^2\phi$ if and only if $\Tr{\ad_{p_i}}=0$ for all $i=1,\ldots,m$.

Fix $g\in G$: due to Lie bracket generating condition, the control system \r{e-controllo} is controllable, then there exists a choice of piecewise constant controls $\fz{u_i}{[0,T]}{\R}$ such that the corresponding solution $\ga(.)$ is an horizontal curve steering $\Id$ to $g$. Then $\dot\ga$ is defined for all $t\in[0,T]$ except for a finite set $E$ of switching times.

Consider now the modular function along $\ga$, i.e. $\fmod(\ga(t))$, that is a continuous function, differentiable for all $t\in[0,T]\backslash E$. We compute its derivative using \r{e-derAd}:  we have $\frac{d}{dt}_{|_{t=t_0}}\det(\Ad_{\ga(t)})=\Tr{\ad_\eta}$ with $\eta=\ga(t_0)^{-1}\dot{\ga}(t_0)$. Due to horizontality of $\ga$, we have $\eta=\sum_{i=1}^m a_i p_i$, hence $\Tr{\ad_\eta}=\sum_{i=1}^m a_i \Tr{\ad_{p_i}}=0$. Then the modular function is piecewise constant along $\ga$. Recalling that it is continuous, we have that it is constant. Varying along all $g\in G$ and recalling that $\fmod(\Id)=1$, we have $\fmod\equiv 1$, hence $G$ is unimodular.\qed

\subsubsection{The hypoelliptic Laplacian on a non-unimodular Lie group}
\newcommand{\Nunim}{A^+(\R)\oplus \R}
In this section we present a non-unimodular Lie group endowed with a left-invariant sub-Riemannian structure. We then compute the explicit expression of the intrinsic hypoelliptic Laplacian: from Proposition \ref{p-inversa} we already know that it is the sum of squares plus a first order term.

Consider the Lie group $$\Nunim:=\Pg{\Pt{\ba{ccc}
a&0&b\\
0&1&c\\
0&0&1
\ea}\ |\ a>0,\,b,c\in\R}.$$
It is the direct sum of the group $A^+(\R)$ of affine transformations on the real line $x\mapsto ax+b$ with $a>0$ and the additive group $(\R,+)$. Indeed, observe that $$\Pt{\ba{ccc}
a&0&b\\
0&1&c\\
0&0&1
\ea}\Pt{\ba{c}x\\d\\1\ea}=\Pt{\ba{c}ax+b\\c+d\\1\ea}.$$
The group is non-unimodular, indeed a direct computation gives $\mu_L=\frac{1}{a^2}da\,db\,dc$ and $\mu_R=\frac{1}{a}da\,db\,dc$.

Its Lie algebra $a(\R)\oplus\R$ is generated by 
\bqn
p_1=\Pt{\ba{ccc}
1&0&0\\
0&0&0\\
0&0&0
\ea},\qquad p_2=\Pt{\ba{ccc}
0&0&1\\
0&0&1\\
0&0&0
\ea},\qquad k=\Pt{\ba{ccc}
0&0&1\\
0&0&0\\
0&0&0
\ea},
\eqnn
for which the following commutation rules hold: $[p_1,p_2]=k\ [p_2,k]=0\ [k,p_1]=-k$.

We define a trivializable sub-Riemannian structure on $\Nunim$ as presented in Section \ref{ss-leftmanifold}: consider the two left-invariant vector fields $X_i(g)=g p_i$ with $i=1,2$ and define
\bqn
\distr(g)=\span{X_1(g),X_2(g)}~~~~~~\metr_g(X_i(g),X_j(g))=\de_{ij}.
\eqnn
Using \r{e-Dh-Lie}, one gets the following expression for the hypoelliptic Laplacian:
\bqn
\Dh \phi = L_{X_1}^2\phi + L_{X_2}^2 \phi + L_{X_1} \phi.
\eqnn

\subsection{The intrinsic hypoelliptic Laplacian in terms of the formal adjoints of the vector fields}

In the literature \gr{another common} definition of hypoelliptic Laplacian \gr{can be found} (see for instance \cite{jerison2}):
\bqn
\llabel{LastL}
\Delta^*=-\sum_{i=1}^m L_{X_i}^\ast L_{X_i},
\eqn
where $\{X_1,\ldots, X_m\}$ is a set of vector fields satysfying the H\"ormander condition  and  the formal adjoint $L_{X_i}^\ast$ is computed with respect to a given volume form. This expression clearly simplifies to the sum of squares when the vector fields are formally skew-adjoint, i.e. $L_{X_i}^\ast =-L_{X_i}$.

In this section we show that our definition of intrinsic hypoelliptic Laplacian coincides locally with \r{LastL}, when 
 $\{X_1,\ldots, X_m\}$ is an orthonormal frame for the sub-Riemannian manifold and the formal adjoint of the vector fields are computed with respect to the sub-Riemannian volume form. 
 %This fact should be evident for experts, but to have a self-consistent paper, we give the proof.
  
We then show that left-invariant vector fields on a Lie group $G$ are formally skew-adjoint with respect to the right-Haar measure, providing an alternative proof of the fact that for unimodular Lie Groups the intrinsic hypoelliptic Laplacian is the sum of squares.\\\\
\bp Locally, 
%and on functions belonging to $\con^\infty_c(M,\R)$ 
the  intrinsic hypoelliptic Laplacian $\Dh$ can be written as
$-\sum_{i=1}^m L_{X_i}^\ast L_{X_i},$ where $\{X_1,\ldots, X_m \}$ is a local orthonormal frame, and $L_{X_i}^\ast$ is the formal adjoint of the Lie derivative $L_{X_i}$ of the vector field $X_i$, i.e.
\bqn
(\phi_1,L_{X_i}^\ast,\phi_2)=(\phi_2,L_{X_i}\phi_1),\mbox{ for every }\phi_1,\phi_2\in \con^\infty_c(M,\R),~~i=1,\ldots,m,
\eqn
and the scalar product is the one of $L^2(M,\R)$ with respect to the invariant volume form, i.e. $(\phi_1,\phi_2):=\int_M\phi_1 \,\phi_2\,\muh$.
\ep
\proof
Given a volume form $\mu$ on $M$, a definition of divergence of a smooth vector field $X$ (equivalent to $L_X\mu=\div(X)\mu$) is
\bqn
\int_M \div(X) \phi\,\mu=-\int_M L_X \phi \mu, \mbox{ for every }\phi\in \con^\infty_c(M,\R);\nonumber
\eqn 
see for instance \cite{taylor}. 
We are going to prove that 
\bqn
\Dh \phi=-\sum_{i=1}^m L_{X_i}^\ast L_{X_i}\phi,~~
\mbox{ for every }\phi\in\con^\infty_c(M,\R):
\eqnl{unouno}
indeed, multiplying the left-hand side of \r{unouno} by $\psi\in\con_c^\infty(M)$  and integrating with respect to $\muh$ we have,
\bqn
\int_M \left( \Dh\phi \right) \psi\,\muh&=&    
\int_M \left(  \divh (\gradh \phi)\right) \psi\,\muh
= \int_M \divh \left(\sum_{i=1}^n (L_{X_i}\phi) {X_i}\right) \psi \,\muh=\nn
&=&-\int_M  \sum_{i=1}^n (L_{X_i}\phi) (L_{X_i}\psi)\,\muh.
\eqnn
For the right hand side we get the same expression. Since $\psi$ is arbitrary, the conclusion follows.

Then, by density, one concludes that $\Dh \phi=-\sum_{i=1}^m L_{X_i}^\ast L_{X_i}\phi,~~
\mbox{ for every }\phi\in\con^2(M,\R)$.
\qed
\bp
Let $G$ be a Lie group and $X$ a left-invariant vector field on $G$. Then $L_X$ is formally skew-adjoint with respect to the right-Haar measure.
\ep
\proof
Let $\phi\in\con^\infty_c(M,\R)$ and $X=gp$ ($p\in\l,$ $g\in G$). Since $X$ is left-invariant and $\mu_R$ is right-invariant, we have
\bqn
\int_G (L_X \phi)(g_0) \,\mu_R(g_0)&=&
\int_G\left.\frac{d}{dt}\right|_{t=0}\phi(g_0e^{t p})\,\mu_R(g_0)=
\left.\frac{d}{dt}\right|_{t=0}\int_G\phi(g_0e^{t p})\,\mu_R(g_0)\nn
&=&\left.\frac{d}{dt}\right|_{t=0}\int_G\phi(g') \mu_R(g'e^{-tp})=\left.\frac{d}{dt}\right|_{t=0}\int_G\phi(g') \mu_R(g')=0.
\eqnn
Hence, for every $\phi_1,\phi_2\in\con^\infty_c(M,\R)$ we have
$$
 0=\int_G L_X(\phi_1 \phi_2)\,\mu_R= \int_G L_X(\phi_1) \,\phi_2\,\mu_R+ \int_G \phi_1 \,(L_X \phi_2)\,\mu_R=(\phi_2,L_X\phi_1)+(\phi_1,L_X\phi_2)
$$
and the conclusion follows.
\qed\\
For unimodular groups we can assume $\muh=\mu_L=\mu_R$  (cfr. Remark \ref{r-misure}) and left-invariant vector fields are formally skew-adjoint with respect to $\muh$. This argument provides an alternative proof of the fact that on unimodular Lie groups the hypoelliptic Laplacian is the sum of squares.

\section{The Generalized Fourier Transform on unimodular Lie groups}
\label{s-genfour}
\newcommand{\Rep}{\Repbase^\lam}
\newcommand{\dRep}{\dRepbase^\lam}
\newcommand{\domR}{\domRbase^\lam}
\newcommand{\fzR}{\fzRbase^\lam_n}
\newcommand{\eigR}{\eigRbase^\lam_n}

Let $f\in L^1(\R,\R)$: its Fourier transform is defined by the formula
$$
\hat f(\lam)=\int_\R f(x) {e^{-i x \lam}} dx.
$$
If $f\in L^1(\R,\R) \cap L^2(\R,\R)$ then $\hat f\in L^2(\R,\R)$ and one has
$$
\int_\R|f(x)|^2dx=\int_\R|\hat f(\lam)|^2  \frac{d\lam}{2\pi},
$$
called Parseval or Plancherel equation. By density of $L^1(\R,\R) \cap L^2(\R,\R)$ in $L^2(\R,\R)$, this equation expresses the fact that the Fourier transform is an isometry between $L^2(\R,\R)$ and itself. Moreover, the following inversion formula holds:
$$
f(x)=\int_\R \hat f(\lam) e^{i x \lam} \frac{d\lam}{{2 \pi}},
$$
where the equality is intended in the $L^2$ sense. 
It has been known from more than 50 years that the Fourier transform generalizes to a wide class of locally compact groups (see for instance \cite{chirichian,duflo,hewittI,hewittII,kirillov-el,taylor-GFT}).  Next we briefly present this generalization for groups satisfying the following \hp:
\bd
\i  \Hp\ $G$ is a unimodular Lie group of Type I.
\ed
For the definition of groups of Type I see \cite{dixmier1}. For our purposes it is sufficient to recall that all groups treated in this paper (i.e. \gruppi) are of Type I. \gr{Actually,} {\it both the real connected semisimple and the real connected nilpotent Lie groups are of Type I} \cite{dixmier,harish} and even though not all solvable groups are of Type I, this is the case for the group of the rototranslations of the plane  $\mot$ \cite{sugiura}. 
In the following, the $L^p$ spaces $L^p(G,\C)$ are intended with respect to the Haar measure $\mu:=\mu_L=\mu_R$.
 
Let $G$ be a Lie group satisfying \Hp\ and $\hat G$ be the dual\ffoot{See footnote \ref{nota-dual}.} of the group $G$, i.e. the set of all equivalence classes of unitary irreducible representations of $G$. Let $\lam\in\hat G$: in the following we indicate by $\Rep$ a choice of an irreducible representation in the class $\lam$. By definition, $\Rep$ is a map that to an element of $G$ associates a unitary operator acting on a complex separable  Hilbert space ${\H^\lam}$:
\mmfunz{\Rep}{G}{U(\domR)}{g}{\Rep(g).}

The index $\lam$ for $\H^\lam$ indicates that in general the Hilbert space can vary with $\lam$.
 
\bdeff
Let $G$ be a Lie group satisfying \Hp, and $f\in L^1(G,\C)$. The generalized (or noncommutative) Fourier transform (\GFT) of $f$ is the map (indicated in the following as $\hat f$ or $\FF(f)$)
that to each element of $\hat G$ associates the linear operator on $\domR$:
 \bqn
 \hat f(\lam):=\FF(f):=\int_G    f(g)\Rep(g^{-1})d\mu.
 \eqn
  \edeff
\noindent
Notice that since $f$ is integrable and  $\Rep$ unitary, then $\hat f(\lam)$ is a bounded operator.

\brem
$\hat f$ can be seen as an operator from $\stackrel{\oplus}{\int_{\hat G}}
\H^\lam$ to itself. We also use the notation  $\hat f=\stackrel{\oplus}{\int_{\hat G}}\hat f(\lam)$
\erem

In general $\hat G$ is not a group and its structure can be quite complicated.  In the case in which $G$ is abelian then $\hat G$ is a group; if $G$ is nilpotent then $\hat G$ has the structure of $\R^n$ for some $n$; if $G$ is compact then it is a Tannaka category
(moreover, in this case each ${\H}^\lam$ is finite dimensional). Under the \hp\ \Hp\ one can define on $\hat G$ a positive measure $dP(\lam)$   (called the Plancherel measure)  such that for every $f\in L^1(G,\C)\cap L^2(G,\C)$ one has 
$$
\int_G |f(g)|^2 \mu(g)=\int_{\hat G}Tr(\hat f(\lam)\circ \hat f (\lam)^\ast ) dP(\lam).
$$
By density of $L^1(G,\C)\cap L^2(G,\C)$ in $L^2(G,\C)$, this formula expresses the fact that the \GFT\ is an isometry between $L^2(G,\C)$ and 
$\stackrel{\oplus}{\int_{\hat G}} \HS^\lam$, the set  of Hilbert-Schmidt operators with respect to the Plancherel measure.
Moreover, it is obvious that:
\bp
Let $G$ be a Lie group satisfying \Hp and $f\in L^1(G,\C)\cap L^2(G,\C)$. We have, for each $g\in G$
\bqn
f(g)=\int_{\hat G}Tr(\hat f (\lam)\circ   \Rep(g) ) dP(\lam).
\eqn  where the equality is intended in the $L^2$ sense. 
\ep
\noindent
{It is immediate to verify that, given two functions $f_1,f_2\in  L^1(G,\C)$ and defining their convolution as
\bqn
(f_1\ast f_2)(g)=\int_G  f_1(h)f_2(h^{-1}g) dh,
\eqn
then the GFT maps  the convolution into non-commutative product:  
\bqn 
\FF( f_1\ast f_2) (\lam)=\hat{f}_2 (\lam)\hat{f}_1(\lam).
\label{conv}
\eqn
Another important property is that if $\delta_{\Id}(g)$ is the Dirac function at the identity over $G$, then 
\bqn
\hat{\delta}_{\Id}(\lam)=\Id_{H^\lam}.
\eqn
}
In the following, a key role is played by  the differential of the representation $\Rep$, that is the map 
\bqn
\label{diff-rep}
d\Rep:X\mapsto d\Rep(X):=\left.\frac{d}{dt}\right|_{t=0}\Rep(e^{t p}),
\eqn
where $X=gp$, ($p\in\l$, $g\in G$) is a left-invariant vector field over $G$. By Stone theorem (see for instance \cite[p. 6]{taylor-GFT})  $d\Rep(X)$ is a (possibly unbounded) skew-adjoint operator on 
%the tangent of $\H_\lam$ that we identify with 
$\H^\lam$. 
We have the following:
\bp
\label{p-uffa}
Let $G$ be a Lie group satisfying \Hp\ and $X$ be a left-invariant vector field over $G$. The GFT of $X$, i.e. $\hat X=\FF L_X\FF^{-1}$ splits into the Hilbert sum of operators $\hat X^\lam$, \gr{each one of which acts} on the set $\HS^\lam$ of Hilbert-Schmidt operators over $\H^\lam$:
  \bqn
 \hat X=\stackrel{\oplus}{\int_{\hat G}}\hat X^\lam.
 \eqnn
Moreover,
  \bqn
 \hat X^\lam\Xi=d\Rep(X)\circ\Xi,~~\mbox{ for every }\Xi\in \HS^\lam,
 \eqn
 i.e. the GFT of a left-invariant vector field acts as a left-translation over $\HS^\lam$.
\ep
\noindent
\proof
Consider the GFT of the operator $R_{e^{t p}}$ of right-translation of a function by ${e^{t p}}$, $p\in\l$, i.e.
$$
\left(R_{e^{t p}}f\right)(g_0)=f(g_0e^{t p}),
$$
and compute its GFT:
\bqn
\FF\left(R_{e^{t p}}f\right)(\lam)&=&\FF\left(f(g_0e^{t p})\right)(\lam)=
\int_G  f(g_0e^{t p})\Rep(g_0^{-1}) \mu(g_0)
= \int_G    f(g')\Rep(e^{t p})\Rep(g'^{-1}) \mu(g' e^{-t p})\nn
&=&
\left( \Rep(e^{t p})\right) \hat f(\lam),
\nonumber
\eqnn
where in the last equality we use the right-invariance of the Haar measure.
Hence the GFT acts as a left-translation on $\HS^\lam$ and it disintegrates the right-regular representations. It follows that
$$
\hat R_{e^{tp}}=\FF R_{e^{tp}}\FF^{-1}=\stackrel{\oplus}{\int_{\hat G}}\Rep(e^{tp}).
$$
Passing to the infinitesimal generators, with $X=gp$, the conclusion follows. \qed

\brem
\label{r-accalambda}
From the fact that the GFT of a left-invariant vector field acts as a left-translation, it follows that $\hat X^\lam$ can be interpreted as an operator over $\H^\lam$.
\erem

\subsection{Computation of the kernel of the hypoelliptic heat equation}

\newcommand{\du}{L}

In this section we provide a general method  to compute the kernel of the hypoelliptic heat equation on a left-invariant sub-Riemannian manifold $(G,\distr,\metr)$ such that $G$ satisfies the assumption \Hp. 

We begin by recalling some existence results (for the semigroup of evolution and for the corresponding kernel) in the case of the  sum of squares. We recall that for all the examples treated in this paper, the invariant hypoelliptic Laplacian is the sum of squares. 

Let $G$ be a unimodular Lie group and $(G,\distr,\metr)$ a left-invariant sub-Riemannian manifold generated by the orthonormal basis $\{p_1,\ldots,p_m\}$, and consider the hypoelliptic heat equation
\bqn
\partial_t\phi(t,g)=\Dh\phi(t,g).
\llabel{eq-hypoQ}
\eqn
Since $G$ is unimodular, then $\Dh=L_{X_1}^2+\ldots +L_{X_m}^2$, where $L_{X_i}$ is the Lie derivative w.r.t. the vector field $X_i:=gp_i$ $(i=1,\ldots, m)$.
Following Varopoulos \cite[pp. 20-21, 106]{varopoulos}, since $\Dh$ is a sum of squares, then it is a symmetric operator that we identify with its Friedrichs (self-adjoint) extension, that is the infinitesimal generator of a (Markov) semigroup $e^{t \Dh }$. Thanks to the left-invariance of $X_i$ (with $i=1,\ldots, m)$, $e^{t \Dh }$ admits a a right-convolution kernel $p_t(.)$, i.e. 
\bqn
 e^{t \Dh } \phi_0(g)=\phi_0\ast p_t (g)=\int_G  \phi_0(h)p_t(h^{-1}g) \mu(h)
\eqn
is the solution for $t>0$ to $\r{eq-hypoQ}$ with initial condition $\phi(0,g)=\phi_0(g)\in L^1(G,\R)$ with respect to the Haar measure.

Since the operator $\partial_t-\Dh$ is hypoelliptic, then the kernel is a $\con^\infty$ function of $(t,g)\in \R^+\times G$. 
Notice that $p_t(g)=e^{t \Dh }\de_\Id(g)$.

The main results of the paper are based on the following key fact.
\bt
\label{t-main}
Let $G$ be a Lie group satisfying \Hp\ and $(G,\distr,\metr)$  a left-invariant sub-Riemannian manifold generated by the orthonormal basis $\{p_1,\ldots,p_m\}$. Let $\Dh=L_{X_1}^2+\ldots +L_{X_m}^2$ be the intrinsic hypoelliptic Laplacian where $L_{X_i}$ is the Lie derivative w.r.t. the vector field $X_i:=gp_i$. 

Let $\Pg{\Rep}_{\lam\in \hat G}$ be the set of all non-equivalent classes of irreducible representations of the group $G$, each acting on an Hilbert space $\domR$, and $dP(\lam)$ be the Plancherel measure on the dual space $\hat G$. We have the following:
\be[{\bf i)}]
\i the GFT of $\Dh$ splits into the Hilbert sum of operators $\DHL$, \gr{each one of which leaves}  ${\H^\lam}$ invariant:
\bqn
\llabel{oplus}
\hat{\Delta}_\SR=\FF\Dh\FF^{-1}=\stackrel{\oplus}{\int_{\hat G}}\DHL dP(\lam),
\mbox{~~ where ~~}
\DHL=\sum_{i=1}^m \left(\hat{X}_i^\lam\right)^2.
\eqn
\i The operator $\DHL$ is self-adjoint and it is the infinitesimal generator of a contraction semi-group $e^{t\DHL}$ over $\HS^\lam$, i.e. $e^{t\DHL}  \Xi_0^\lam$ is the solution for $t>0$ to the operator equation $\partial_t  \Xi^\lam(t)=\DHL \Xi^\lam(t)$ in $\HS^\lam$, with initial condition  $\Xi^\lam(0)= \Xi^\lam_0 $.
\i The hypoelliptic heat kernel is
\bqn
\llabel{formula-1}
p_t(g)=\int_{\hat G}Tr\left(e^{t\DHL}\Rep(g)\right)dP(\lam),~~t>0.
\eqn
\ee
\et
\noindent
\proof
%Let us compute the GFT of $L_Xf$ where $X$ is a left-invariant vector field and $f$ is a smooth integrable function $f:G\to\R$. We have
%\bqn
%\FF (L_X f)(\lam) &=&\FF\left(   \left.\frac{d}{dt}\right|_{t=0} f(g_0e^{t X}) \right)(\lam)=
%\int_G  \left.\frac{d}{dt}\right|_{t=0} f(g_0e^{t X})\Rep(g_0^{-1}) \mu(g_0)\nn
%&=& \left.\frac{d}{dt}\right|_{t=0}\int_G    f(g')\Rep(e^{t X})\Rep(g'^{-1}) \mu(g' e^{-t X})=
%\left( \left.\frac{d}{dt}\right|_{t=0}\Rep(e^{t X})\right) \hat f(\lam).
%\nonumber
%\eqn
%where in the last equality we used the right-invariance of the Haar measure.
%The map $d\Rep:X\mapsto \left.\frac{d}{dt}\right|_{t=0}\Rep(e^{t X})$ is an operator on $U(\H_\lam)$, and it is called the differential of the representation. 
%It follows $\hat X=d\Rep (X)$. In other words, since 
%the GFT desintegrates the right-regular representations and the left-invariant vector fields $X_i$ have one parameter group of right-translations $e^{t X_i}$, it follows that $\hat{X}=\FF X \FF^{-1}=
%\stackrel{\oplus}{\int} X dP(\lam).$
%It means that $\hat X(\lam)$ is an operator over $\H_\lam$, and so
%$\DHL$ and  $e^{t\DHL}$ are as well.
Following Varopoulos as above, and using Proposition \ref{p-uffa}, {\bf i)} follows. 
Item {\bf ii)} follows from the split \r{oplus} and from the fact that  GFT is an isometry between $L^2(G, \C)$ (the set of square integrable functions from $G$ to $\C$ with respect to the Haar measure) and the set  $\stackrel{\oplus}{\int_{\hat G}} \HS^\lam$ of   Hilbert-Schmidt operators with respect to the Plancherel measure.
Item {\bf  iii)} is obtained applying the inverse GFT to $e^{t\DHL} \Xi_0^\lam $ and the convolution formula \r{conv}.
The integral is convergent by the existence theorem for $p_t$, see  \cite[p. 106]{varopoulos}.
\qed
\brem
As a consequence of Remark \ref{r-accalambda}, it follows that $\DHL$ and $e^{t\DHL}$ can be considered as operators on $\H^\lam$.

\erem

%\brem
%\label{ispazi}
%Notice that $\DHL$ is an operator on $U(\H_\lam)$ i.e.  $\DHL\hat f(\lam)$ (where $f\in L^1(G,\C)$)) is an operator on $\H_\lam$. Hence  $\DHL $ can be though as a self-adjoint operator acting on $\H_\lam$ on vectors of the type  $\Xi^\lam=\hat f(\lam)\psi^\lam$, where $\psi^\lam\in\H_\lam$. The same remark holds for $e^{t\DHL}$.
%In the following, when we speak of eigenvalues and eigenvectors of $\DHL$, we think to  $\DHL$ as an operator on  $\H_\lam$. 
%\erem

\gr{In the case when each $\DHL$ has discrete spectrum, the following corollary gives an explicit formula for the hypoelliptic heat kernel in terms of its eigenvalues and eigenvectors.}

\bc Under the hypotheses of Theorem \ref{t-main}, if in addition \gr{we have that} for every $\lam$, $\DHL$ (considered as an operator over $\H^\lam$) has discrete spectrum and  $\Pg{\fzR}$ is a complete set of eigenfunctions of norm one with the corresponding set of eigenvalues $\Pg{\eigR}$, then
\bqn
p_t(g)=\int_{\hat G}\left( \sum_{n} e^{\eigR t} \Pa{\fzR ,\Rep(g) \fzR} \right)d P(\lam)
\eqnl{e-fundsol-general}
where $\Pa{.,.}$ is the scalar product in $\domR$.
\ec
\proof
Recall that $\Tr{AB}=\Tr{BA}$ and that $\Tr{A}=\sum_{i\in I} \Pa{e_i,A e_i}$ for any complete set $\Pg{e_i}_{i\in I}$ of orthonormal vectors. Hence $\Tr{e^{t\DHL}\Rep(g)}=\sum_{n} \Pa{\fzR ,\Rep(g) e^{t\DHL} \fzR}$. Observe that $\partial_t \fzR=\DHL \fzR=\eigR \fzR$, hence $e^{t\DHL} \fzR=e^{\eigR t}\fzR$, from which the result follows.
\qed \\\\

The following corollary gives a useful formula for the hypoelliptic heat kernel in the case in which for all $\lam\in\hat{G}$ each operator $e^{t\DHL}$ admits a convolution kernel $Q_t^\lam(.,.)$. Here by $\fzRbase^\lam$, we \gr{intend} an element of $\H^\lam$.

\bc
\label{c-2} Under the hypotheses of Theorem \ref{t-main}, if for all $\lam\in\hat G$ we have $\H^\lam=L^2(X^\lam,d\th^\lam)$ for some measure space $(X^\lam,d\th^\lam)$ and 
$$\Pq{e^{t\DHL} \fzRbase^\lam}(\th)=\int_{X^\lam}  \fzRbase^\lam (\bar \th) Q_t^\lam(\th,\bar\th)\,d\bar \th,$$ then
%%%%%%%%%%%%%%%%%%%%%%%%%%%
$$
p_t(g)=\int_{\hat G}\int_{X^\lam}\left.\Rep(g)Q_t^\lam(\th,\bar\th)\right|_{\th=\bar \th}              \,d\bar \th\,
d P(\lam),
$$
where in the last formula $\Rep(g)$ acts on $Q_t^\lam(\th,\bar\th)$ as a function of $\th$.
\ec
\proof  
From \r{formula-1}, we have
$$
p_t(g)=\int_{\hat G}Tr\left(e^{t\DHL}\Rep(g)\right)dP(\lam)=\int_{\hat G}Tr\left(\Rep(g)e^{t\DHL}\right)dP(\lam).
$$
We have to compute the trace of the operator
\newcommand{\oprusso}{\Theta}
\bqn
\label{kernellino}
\oprusso=\Rep(g)e^{t\DHL}:\fzRbase^\lam(\th)\mapsto \Rep(g)e^{t\DHL}\fzRbase^\lam(\th)&=&
\Rep(g)\int_{X^\lam}  \fzRbase^\lam(\bar \th) Q_t^\lam(\th,\bar\th)\,d\bar \th=\\
&=&\int_{X^\lam}  K(\th,\bar\th) \fzRbase^\lam(\bar \th)d\bar \th
\eqnn
where  $K(\th,\bar\th)= \Rep(g)Q_t^\lam(\th,\bar\th)$ is a function of $\th,\bar\th$ and $\Rep(g)$ acts on $Q_t^\lam(\th,\bar\th)$ as a function of $\th$. The trace of $\oprusso$ is $\int_X K(\bar \th,\bar \th)d\bar\th$ and the conclusion follows.
\qed

\section{Explicit expressions on 3D unimodular Lie groups}
\label{s-examples}
\renewcommand{\dRepbase}[2]{\hat{X}_{#1}^{#2}}
\newcommand{\DHLbase}[1]{\hat{\Delta}_\SR^{#1}}

\subsection{The hypoelliptic heat equation on $H_2$}

\renewcommand{\Rep}{\Repbase^\lam}
\renewcommand{\dRep}[1]{\dRepbase{#1}{\lam}}
\renewcommand{\domR}{\domRbase}
\renewcommand{\fzR}{\fzRbase}
\renewcommand{\DHL}{\DHLbase{\lam}}

In this section we apply the method presented above to solve the hypoelliptic heat equation \r{eq-hypoQ} on the Heisenberg group. This kernel, via the GFT, was first  obtained  by Hulanicki (see \cite{hul}). We present it as an application of Corollary \ref{c-2}, since in this case an expression for the kernel of the GFT of this equation is known.

We write the Heisenberg group as the 3D group of matrices
$$H_2=\Pg{\Pt{\ba{ccc}
1&x&z+\frac12 x y\\
0&1&y\\
0&0&1
\ea}\ |\ x,y,z\in\R}$$
endowed with the standard matrix product. $H_2$ is indeed $\R^3,$  
$$(x,y,z)\sim\Pt{\ba{ccc}
1&x&z+\frac12 x y\\
0&1&y\\
0&0&1
\ea},$$ endowed with  the group law  $$(x_1,y_1,z_1)\cdot(x_2,y_2,z_2)=\Pt{x_1+x_2,y_1+y_2,z_1+z_2+\frac12\Pt{x_1y_2-x_2y_1}}.
$$
A basis of its Lie algebra  is $\Pg{p_1,p_2,k}$ where
\bqn
p_1=\Pt{\ba{ccc}
0&1&0\\
0&0&0\\
0&0&0
\ea}\quad
p_2=\Pt{\ba{ccc}
0&0&0\\
0&0&1\\
0&0&0
\ea}\quad
k=\Pt{\ba{ccc}
0&0&1\\
0&0&0\\
0&0&0
\ea}.
\eqnl{e-H2-alg}
They satisfy the following commutation rules: $[p_1,p_2]=k$, $[p_1,k]=[p_2,k]=0$, hence $H_2$ is a 2-step nilpotent group. We define a left-invariant sub-Riemannian structure on $H_2$ as presented in Section \ref{ss-leftmanifold}: consider the two left-invariant vector fields $X_i(g)=g p_i$ with $i=1,2$ and define
\bqn
\distr(g)=\span{X_1(g),X_2(g)}~~~~~~\metr_g(X_i(g),X_j(g))=\de_{ij}.
\eqnn
Writing the group $H_2$ in coordinates $(x,y,z)$ on $\R^3$, we have the following expression for the Lie derivatives of $X_1$ and $X_2$: 
$$L_{X_1}=\partial_x-\frac y2 \partial_z,\qquad L_{X_2}=\partial_y+\frac x2 \partial_z.$$
The Heisenberg group is unimodular, hence the hypoelliptic Laplacian $\Dh$ is the sum of squares:
\bqn
\Dh\phi=\Pt{L_{X_1}^2+ L_{X_2}^2}\phi.
\label{eq-H2-heat}
\eqn
\brem
It is interesting to notice that all left-invariant sub-Riemannian structures that one can define on the Heisenberg group are isometric.
\erem

In the next proposition we present the structure of the dual group of $H_2$. For details and proofs see for instance \cite{kirillov}.

\bp
The dual space of $H_2$ is $\hat{G}=\Pg{\Rep\ |\ \lam\in\R}$, where
\mmfunz{\Rep(g)}{\domR}{\domR}{\fzR(\th)}{e^{i\lam\Pt{z-y\th+\frac{xy}{2}}}\fzR(\th-x),}
whose domain is $\domR=L^2(\R,\C)$, endowed with the standard product $<\fzR_1,\fzR_2>:=\int_\R \fzR_1(\th)\overline{\fzR_2}(\th)\,d\th$ where $d\th$ is the Lebesgue measure.

The Plancherel measure on $\hat{G}$ is $dP(\lam)=\frac{|\lam|}{4\pi^2} d\lam$, where $d\lam$ is the Lebesgue measure on $\R$.
\ep
\brem
Notice that in this example the domain $\domR$ of the representation $\Rep$ does not depend on $\lam$.
\erem

\subsubsection{The kernel of the hypoelliptic heat equation}

Consider the representation $\Rep$ of $H_2$ and let $\dRep{i}$ be the corresponding representations of the differential operators $L_{X_i}$ with $i=1,2$. Recall that $\dRep{i}$  are operators on $\H$.
%As in Remark \ref{ispazi}, recall that $\dRep{i}$  can be though
%as an operator on $\H$ on functions of the type  $\Xi^\lam(\th)=\hat f(\lam)\psi^\lam(\th)$, where $\psi^\lam\in\H$.  
From formula \r{diff-rep} we have $$[\dRep{1} \fzRbase](\th)=-\frac{d}{d\th}\fzRbase(\th),\quad
[\dRep{2}\fzRbase](\th)=
-i\lam\th \,\fzRbase(\th),\mbox{~~~hence~~}[\DHL\fzRbase](\th)=\Pt{\frac{d^2}{d\th^2}-\lam^2\th^2}\fzRbase(\th).$$
The GFT of the hypoelliptic heat equation is thus
$$\partial_t\fzRbase= \Pt{\frac{d^2}{d\th^2}-\lam^2\th^2}\fzRbase.$$

The kernel of this equation is known (see for instance \cite{beals}) and it is called the Mehler kernel (its computation is very similar to the computation of the kernel for the harmonic oscillator in quantum mechanics):
%One have to find eigenvalues and eigenfunctions of the operator $\frac{d^2}{d\th^2}-\lam^2\th^2$ an using the Mehler formula gets:
$$
Q_t^\lam(\th,\bar{\th}):=\sqrt{\frac{\lam}{2\pi\sinh(2\lam t)}} \exp\left({-\frac12\frac{\lam\cosh(2\lam t)}{\sinh(2\lam t)} (\th^2+\bar{\th}^2)+\frac{\lam\th\bar{\th}}{\sinh(2\lam t)}}\right).
$$
Using Corollary \ref{c-2} and after straightforward computations, one gets  the kernel of the hypoelliptic heat equation on the Heisenberg group:
\bqn
p_t(x,y,z)=\frac{1}{(2\pi t)^2}\int_\R \frac{2 \tau}{\sinh(2\tau)}\exp\left({-\frac{\tau(x^2+y^2)}{2t\tanh(2\tau)}}\right)\cos(2\frac{z\tau}{t}) d\tau.
\eqnl{eq-H2-heat-exp}
This formula differs from the one by Gaveau \cite{gaveau} for some numerical factors since he studies the equation 
$$\partial_t\phi=\frac12\left( (\partial_x  +2y\partial_z)^2+(\partial_y  -2x\partial_z)^2 \right)\phi.$$ The Gaveau formula is recovered from \r{eq-H2-heat-exp} with $t\to t/2$ and $z\to z/4$. Moreover, a multiplicative factor $\frac14$ should be added, because from the change of variables one gets that the Haar measure is $\frac14 dx\,dy\,dz$ instead of $dx\,dy\,dz$ as used by Gaveau.

\subsection{The hypoelliptic heat equation on $SU(2)$}
\newcommand{\matSU}{\Mat{cc}{
\alpha&\beta\\
-\bar \beta& \bar \alpha}}

\renewcommand{\Rep}{\Repbase^{n}}
\renewcommand{\dRep}[1]{\dRepbase{#1}{n}}
\renewcommand{\domR}{\domRbase^{n}}
\renewcommand{\fzR}[1][]{\fzRbase^n_{k #1}}
\renewcommand{\DHL}{\DHLbase{n}}

In this section we solve the hypoelliptic heat equation \r{eq-hypoQ} on the Lie group
$$SU(2)=\Pg{\matSU\ |\  \al,\beta\in\C, |\alpha|^2+|\beta|^2=1}.$$

A basis of the Lie algebra $su(2)$ is $\Pg{p_1,p_2,k}$ where\ffoot{See \cite[pp. 67]{sugiura}: $p_1=X_1$, $p_2=X_2$, $k=X_3$.}
\bqn
p_1=\frac{1}{2}\Pt{\begin{array}{cc}
0 & i\\
i & 0
\end{array}}\quad
p_2=\frac{1}{2}\Pt{\begin{array}{cc}
0 & -1\\
1 & 0
\end{array}}\quad
k=\frac{1}{2}\Pt{\begin{array}{cc}
i & 0\\
0 & -i
\end{array}}.
\eqnl{pauli}

We define a trivializable sub-Riemannian structure on $SU(2)$ as presented in Section \ref{ss-leftmanifold}: consider the two left-invariant vector fields $X_i(g)=g p_i$ with $i=1,2$ and define
\bqn
\distr(g)=\span{X_1(g),X_2(g)}~~~~~~\metr_g(X_i(g),X_j(g))=\de_{ij}.
\eqnn

The group $SU(2)$ is unimodular, hence the hypoelliptic Laplacian $\Dh$ has the following expression:
\bqn
\Dh\fzRbase=\Pt{L_{X_1}^2+ L_{X_2}^2}\fzRbase
\label{eq-SU2-heat}
\eqn

In the next proposition we present the structure of the dual group of $SU(2)$. For details and proofs see for instance \cite{sugiura}.

\bp
The dual space of $SU(2)$ is $\hat{G}=\Pg{\Rep\ |\ n\in\N}$.

The domain $\domR$ of $\Rep$ is the space of homogeneous polynomials of degree $n$ in two variables $(z_1,z_2)$ with complex coefficients $\domR:=\Pg{\sum_{k=0}^{n}a_kz_1^kz_2^{n-k}\ |\ a_k\in\C}$, endowed with the scalar product $$\Pa{\sum_{k=0}^{n}a_kz_1^kz_2^{n-k},\sum_{k=0}^{n}b_kz_1^kz_2^{n-k}}:=\sum_{k=0}^{n}\ k!\ (n-k)!\ a_k\  \bar{b}_k.$$

The representation $\Rep$ is
\mmfunz{\Rep(g)}{\domR}{\domR}{\sum_{k=0}^{n}a_kz_1^kz_2^{n-k}}{\sum_{k=0}^{n}a_kw_1^kw_2^{n-k}}
with $(w_1,w_2)=(z_1,z_2)g=(\al z_1-\bar{\beta}z_2, \beta z_1 +\bar{\al} z_2)$.

The Plancherel measure on $\hat{G}$ is $dP(n)=(n+1)\cont(n)$, where $\cont$ is the counting measure.
\ep

Notice that an orthonormal basis of $\domR$ is $\Pg{\fzR}_{k=0}^n$ with $\fzR:=\frac{z_1^kz_2^{n-k}}{\sqrt{k!\ (n-k)!}}$.

\subsubsection{The kernel of the hypoelliptic heat equation}

Consider the representations $\dRep{i}$ of the differential operators $L_{X_i}$  with $i=1,2$: they are operators on $\domR$, whose action on the basis $\Pg{\fzR}_{k=0}^n$ of $\domR$ is (using formula \r{diff-rep}) 
$$\dRep{1}\fzR=\frac{i}{2}\Pg{k\fzR[-1]+(n-k)\fzR[+1]}\qquad
\dRep{2}\fzR=\frac{1}{2}\Pg{k\fzR[-1]-(n-k)\fzR[+1]}$$
hence $\DHL\fzR=\Pt{k^2-kn-\frac{n}{2}}\fzR.$ Thus, the basis $\Pg{\fzR}_{k=0}^n$ is a complete set of eigenfunctions of norm one for the operator $\DHL$. We are now able to compute the kernel of the hypoelliptic heat equation using formula \r{e-fundsol-general}. 

\bp
The kernel of the hypoelliptic heat equation on $(SU(2),\distr,\metr)$ is
\bqn
\llabel{eq-SU2-heat-exp}
p_t(g)=\sum_{n=0}^\infty (n+1) \sum_{k=0}^n 
e^{(k^2-kn-\frac{n}{2})t} A^{n,k}(g),
\eqn
where $$A^{n,k}\Pt{g}:=\Pa{\fzR,\Rep(g)\fzR}= \sum_{l=0}^{\min\Pg{k,n-k}}
\Pt{\ba{c} k\\k-l\ea}
\Pt{\ba{c} n-k\\l\ea}
\bar{\al}^{k-l} {\al}^{n-k-l}\Pt{|\alpha|^2 -1 }^l$$ with $g=\Pt{\ba{cc}\al&\beta\\-\bar{\beta}&\bar{\al}\ea}$.
\ep
\proof The formula $p_t(g)=\sum_{n=0}^\infty (n+1) \sum_{k=0}^n 
e^{(k^2-kn-\frac{n}{2})t}\Pa{\fzR,\Rep(g)\fzR}$ is given by applying formula \r{e-fundsol-general} in the $SU(2)$ case.

We now prove the explicit expression for $\Pa{\fzR,\Rep(g)\fzR}$: a direct computation gives
\bqn\Rep(g)\fzR
&=&\frac{\sum_{s=0}^{n} \fzRbase_s^n\sqrt{s!\Pt{n-s}!}\Pt{\sum_{l=\max\Pg{0,s-k}}^{\min\Pg{s,n-k}}\Pt{\ba{c}k\\s-l\ea}\Pt{\ba{c}n-k\\l\ea}
\al^{s-l}(-\bar{\beta})^{k-s+l}\beta^l\bar{\al}^{n-k-l}}}{\sqrt{k!\Pt{n-k}!}}
\eqnn
Observe that $\fzR$ is an orthonormal frame for the inner product: hence$$\Pa{\fzR,\Rep(g)\fzR}=\Pa{\fzR,\fzR
\sum_{l=0}^{\min\Pg{k,n-k}}\Pt{\ba{c}k\\k-l\ea}\Pt{\ba{c}n-k\\l\ea}
\al^{k-l}(-\bar{\beta})^{l}\beta^l\bar{\al}^{n-k-l}}.$$ The result easily follows.
\qed

\brem
Notice that, as the sub-Riemannian distance (computed in \cite{nostro-gruppi}), $p_t(g)$ does not 
depend on $\beta$. This is due to the cylindrical symmetry of the distribution around $e^\k=\Pg{e^{ck}\ |\ c\in\R}$.
\erem

\subsection{The hypoelliptic heat equation on $SO(3)$}

\renewcommand{\Rep}{\Repbase^r}
\renewcommand{\dRep}[1]{\dRepbase{#1}{r}}
\renewcommand{\domR}{\domRbase^r}
\renewcommand{\fzR}{\phi^r_s}
\renewcommand{\eigR}{\eigRbase^r_s}
\newcommand{\elem}{g}

Let $\elem$ be an element of $SO(3)=\Pg{A\in\mathrm{Mat}(\R,3)\ |\ AA^T=\Id,~~~\det(A)=1}$. A basis of the Lie algebra $so(3)$ is $\Pg{p_1,p_2,k}$ where\footnote{See \cite[pp. 88]{sugiura}: $p_1=Z_1$, $p_2=Z_2$, $k=Z_3$.}
\bqn
p_1=\Pt{\ba{ccc}
0 & 0 & 0\\
0 & 0 & -1\\
0 & 1 & 0
\ea}\quad
p_2=\Pt{\ba{ccc}
0 & 0 & 1\\
0 & 0 & 0\\
-1 & 0 & 0
\ea}\quad
k=\Pt{\ba{ccc}
0 & -1 & 0\\
1 & 0 & 0\\
0 & 0 & 0
\ea}
\eqnl{eq-baseso3}

We define a trivializable sub-Riemannian structure on $SO(3)$ as presented in Section \ref{ss-leftmanifold}: consider the two left-invariant vector fields $X_i(\elem)=\elem p_i$ with $i=1,2$ and define
\bqn
\distr(\elem)=\span{X_1(\elem),X_2(\elem)}~~~~~~\metr_\elem(X_i(\elem),X_j(\elem))=\de_{ij}.
\eqnn

The group $SO(3)$ is unimodular, hence the hypoelliptic Laplacian $\Dh$ has the following expression:
\bqn
\Dh\phi=\Pt{L_{X_1}^2+ L_{X_2}^2}\phi.
\label{eq-SO3-heat}
\eqn

We present now the structure of the dual group of $SO(3)$. For details and proofs see \cite{sugiura}.

First consider the domain $\domR$, that is the space of complex-valued polynomials of $r$-th degree in three real variables $x,y,z$ that are homogeneous and harmonic
\bqn
\domR=\Pg{f(x,y,z)\ |\ \deg\Pt{f}=r,\ f\mbox{~homogeneous~}, \Delta f= 0}.
\eqnn
Notice that an homogeneous polynomial $f\in\domR$ is uniquely determined by its value on $S^2=\Pg{(x,y,z)\ |\ x^2+y^2+z^2=1}$, as $f(\rho x,\rho y,\rho z)=\rho^r f(x,y,z)$.

Define $\tilde{f}(\al,\beta):=f(\sin(\al)\cos(\beta),\sin(\al)\sin(\beta),\cos(\al))$. Then endow $\domR$ with the scalar product $$<f_1(x,y,z),f_2(x,y,z)>:=\frac{1}{4\pi}\int_{S^2}\tilde{f}_1(\al,\beta)\overline{\tilde{f}_2(\al,\beta)}\sin\al \ d\al\ d\beta.$$

In the following proposition we present the structure of the dual group. 
\bp
The dual space of $SO(3)$ is $\hat{G}=\Pg{\Rep\ |\ r\in\N}$.

Given $\elem\in SO(3)$, the unitary representation $\Rep(\elem)$ is
\mmfunz{\Rep(\elem)}{\domR}{\domR}{f(x,y,z)}{f(x_1,y_1,z_1)}
with $(x_1,y_1,z_1)=(x,y,z)\elem$.

The Plancherel measure on $SO(3)$ is $dP(r)=(2r+1)\cont(r)$, where $\cont$ is the counting measure.
\ep

An orthonormal basis for $\domR$ is given by $\Pg{\fzR}_{s=-r}^r$ with $\tilde{\phi}^r_s(\al,\beta):=e^{is\beta}P_r^{-s}(\cos(\al))$, where $P^s_r(x)$ are the Legendre polynomials\footnote{Recall that $P^s_r(x)$ is defined by $P^s_r(x):=\frac{(1-x^2)^{\frac{s}{2}}}{r! 2^r} \frac{d^{r+s}(x^2-1)^r}{dx^{r+s}}$.}.

\subsubsection{The kernel of the hypoelliptic heat equation}

Consider the representations $\dRep{i}$ of the differential operators $L_{X_i}$  with $i=1,2$: using formula \r{diff-rep} we find the following expressions in spherical coordinates\footnote{i.e. $x=\rho \sin(\al)\cos(\beta)$, $y=\rho \sin(\al)\sin(\beta)$, $z=\rho \cos(\al)$.} 
\bqn
\dRep{1} \fzRbase=
\sin(\beta)\frac{\partial \fzRbase} {\partial\al}+\cot(\al)\cos(\beta)\frac{\partial \fzRbase} {\partial\beta},\qquad \dRep{2} \fzRbase=
-\cos(\beta)\frac{\partial \fzRbase} {\partial\al}+\cot(\al)\sin(\beta)\frac{\partial \fzRbase} {\partial\beta}
\eqnn hence \bqn
\DHL \fzRbase&=&\frac{\partial^2\fzRbase}{\partial\al^2}+\cot^2(\al)\frac{\partial^2\fzRbase}{\partial\beta^2}+
\cot(\al)\frac{\partial\fzRbase}{\partial\al}
\eqn
and its action on the basis $\Pg{\fzR}_{s=-r}^r$ of $\domR$ is
\bqn
\DHL\fzR=\Pt{s^2-r(r+1)} \fzR.
\eqnl{eq-SO3-azioneDelta}
Hence the basis $\Pg{\fzR}_{s=-r}^r$ is a complete set of eigenfunctions of norm one for the operator $\DHL$.
%\proof For the computation of $\dRep{i}$ see \cite{sugiura}. The computation of $\DHL$ is a straightforward consequence.

%For the computation of $\DHL\fzR$, observe that $\DHL=\Delta_{S^2}-\frac{\partial^2}{\partial\beta^2}$, where $\Delta_{S^2}$ is the Laplace-Beltrami operator on $S^2$. Recall that $\Delta_{S^2} \fzR=-n(n+1) \fzR$ as a well known property of spherical harmonics. We also have $\frac{\partial}{\partial\phi} \fzR= im \fzR$, from which formula \r{eq-SO3-azioneDelta} follows. \qed

We compute the kernel of the hypoelliptic heat equation, using \r{e-fundsol-general}.
\bp
The kernel of the hypoelliptic heat equation on $(SO(3),\distr,\metr)$ is
\bqn
p_t(\elem)=\sum_{r=0}^\infty (2r+1) \sum_{s=-r}^r
e^{(s^2-r(r+1))t} <\fzR,\Rep(\elem)\fzR>.
\eqnl{eq-SO3-heat-exp}
\ep

\subsubsection{The heat kernel on $SO(3)$ via the heat kernel on $SU(2)$}

\newcommand{\ul}[1]{\underline{#1}}
\renewcommand{\elem}{\ul{g}}
In this section we verify that the heat kernel on $SO(3)$ given in  \r{eq-SO3-heat-exp} can be easily retrieved from the one on $SU(2)$ given in  \r{eq-SU2-heat-exp}. In the following, all the objects relative to $SO(3)$ are underlined, e.g. $\ul{g}\in SO(3), \ul{p_i}\in so(3),$ the representations $\ul{\Rep}$ acting on $\ul{\domR}$ with basis $\ul{\fzR}$.

Consider the isomorphism of Lie algebras $\fz{\ad}{su(2)}{so(3)}$ defined by $\ad p_1=\ul{p_1},\, \ad p_2=\ul{p_2},\,\ad k=\ul{k}$: it gives the matrix expression of the adjoint map on $su(2)$ with respect to the basis $\Pg{p_1,p_2,k}.$ There is a corresponding endomorphism of groups $\fz{\Ad}{SU(2)}{SO(3)}$ given by $\Ad(\exp v)=\ul{\exp(\ad(v))}$. It is a covering map of $SO(3)$ by $SU(2)$, such that for each matrix $\elem\in SO(3)$ the preimage is given by two opposite matrices $g,-g\in SU(2)$.
\bp
The following relation holds between the kernel $\ul{p_t}$ on $SO(3)$ given in  \r{eq-SO3-heat-exp} and the  kernel $p_t$ on $SU(2)$ given in  \r{eq-SU2-heat-exp}:

$\forall\,\elem\in SO(3), g\in \Ad^{-1}(\elem)$
$$\ul{p_t}(\elem)=\frac{p_t(g)+p_t(-g)}{2}.$$
\ep
\proof
Observe the following key facts (see e.g. \cite[II.\S7]{sugiura}):
\bi
\i on $SU(2)$: $\Repbase^n(-g)\phi=(-1)^n\Repbase^n(g)\phi$;
\i the representation $\ul{\Rep}$ of $SO(3)$ and the representation $\Repbase^{2r}$ of $SU(2)$ are unitarily related, i.e. the following relation holds: $\forall\, g\in SU(2), \ul{\fzR}\in \ul{\domR}$
\bqn
 T^r \ul{\Rep}(\Ad g)\, \ul{\fzR} = \Repbase^{2r} (g) [T^r (\ul{\fzR})].
\eqnl{eq-relSU2SO3}
where the map \mmfunz{T^r}{\ul{\domR}}{\domRbase^{2r}}{\ul{\fzR}}{\phi^{2r}_{r+s}.}
is a unitary isomorphism.
\ei

Then we have explicitly
\bqn
\frac{p_t(g)+p_t(-g)}{2}&=&\frac{\sum_{n=0}^\infty (n+1) (1+(-1)^n)\sum_{k=0}^n 
e^{(k^2-kn-\frac{n}{2})t} A^{n,k}(g)}{2}=\nn
&=&\sum_{r=0}^\infty (2r+1) \sum_{s=-r}^r
e^{(s^2-r(r+1))t} <\phi^{2r}_{r+s},\Repbase^{2r} (g)\phi^{2r}_{r+s}>
\eqnn
after the substitution $r=2n,\ s=k-r$. Using Equation \r{eq-relSU2SO3}, we have $<\phi^{2r}_{r+s},\Repbase^{2r} (g)\phi^{2r}_{r+s}>=<\ul{\fzR},\ul{\Rep}(\Ad g)\ul{\fzR}>$, from which the result directly follows.
\qed

\subsection{The hypoelliptic heat equation on $SL(2)$}
\newcommand{\RepC}[1][j]{\Repbase^{#1,s}}
\newcommand{\fzRC}[1][]{\fzRbase_{m#1}}
\newcommand{\domRC}{\domRbase_C}
\newcommand{\dRepC}[1]{\dRepbase{#1}{j,s}}
\newcommand{\DHLC}{\DHLbase{j,s}}

\newcommand{\RepD}{\Repbase^n}
\newcommand{\fzRD}[1][]{\fzRbase^n_{m#1}}
\newcommand{\domRD}{\domRbase^n}
\newcommand{\dRepD}[1]{\dRepbase{#1}{n}}
\newcommand{\DHLD}{\DHLbase{n}}

\newcommand{\gsu}{\mathscr{G}}

In this section we solve the hypoelliptic heat equation \r{eq-hypoQ} on the Lie group $$SL(2)=\Pg{g\in\mathrm{Mat}(\R,2)\ |\ \det(g)=1}.$$ A basis of the Lie algebra $sl(2)$ is 
$$p_1=\frac{1}{2}\Pt{\begin{array}{cc}
1 & 0\\
0 & -1
\end{array}}\quad
p_2=\frac{1}{2}\Pt{\begin{array}{cc}
0 & 1\\
1 & 0
\end{array}}\quad
k=\frac{1}{2}\Pt{\begin{array}{cc}
0 & -1\\
1 & 0
\end{array}}.$$

We define a trivializable sub-Riemannian structure on $SL(2)$ as presented in Section \ref{ss-leftmanifold}: consider the two left-invariant vector fields $X_i(g)=g p_i$ with $i=1,2$ and define
\bqn
\distr(g)=\span{X_1(g),X_2(g)}~~~~~~\metr_g(X_i(g),X_j(g))=\de_{ij}.
\eqnn

The group $SL(2)$ is unimodular, hence the hypoelliptic Laplacian $\Dh$ has the following expression:
\bqn
\Dh\phi=\Pt{L_{X_1}^2+ L_{X_2}^2}\phi.
\label{eq-SL2-heat}
\eqn

It is well known that $SL(2)$ and $SU(1,1)= \Pg{\Mat{cc}{\al&\beta\\\bar{\beta}&\bar{\al}}\ |\ |\al|^2-|\beta|^2=1}$ are isomorph Lie groups via the isomorphism
$$
\funz{\Pi}{$SL(2)$}{$SU(1,1)$}{$g$}{$\gsu=CgC^{-1}$}
\mbox{~~~~~~~~with~~~} C=\frac{1}{\sqrt{2}}\Mat{cc}{1&-i\\1&i}.$$ This isomorphism also induce an isomorphism of Lie algebras $\fz{d\Pi}{sl(2)}{su(1,1)}$ defined by $d\Pi(p_1)=p_1'$, $d\Pi(p_2)=p_2'$, $d\Pi(k)=k'$ with
$$p_1'=\frac{1}{2}\Pt{\begin{array}{cc}
0 & 1\\
1 & 0
\end{array}}\quad
p_2'=\frac{1}{2}\Pt{\begin{array}{cc}
0 & -i\\
i & 0
\end{array}}\quad
k'=\frac{1}{2}\Pt{\begin{array}{cc}
-i & 0\\
0 & i
\end{array}}.$$
This isomorphism induces naturally the definitions of left-invariant sub-Riemannian structure and of the hypoelliptic Laplacian on $SU(1,1)$.

We present here the structure of the dual of the group $SU(1,1)$, observing that the isomorphism of groups induces an isomorphism of representations. For details and proofs, see \cite{sugiura}.

The dual space $\hat G$ of $SU(1,1)$ contains two continuous and two discrete parts: $\hat G=\hat G_C\disj\hat G_D$ with $\hat G_C=\Pg{\RepC\ |\ j\in\Pg{0,\frac12},s=\frac12+iv, v\in\R^+}$ and $\hat G_D=\Pg{\RepD\ |\ n\in\frac12\Z, |n|\geq 1}$.

We define the domain $\domRC$ of the continuous representation $\RepC$: it is the Hilbert space of $L^2$ complex-valued functions on the unitary circle $S^1=\Pg{x\in\C\ |\ |x|=1}$ with respect to the normalized Lebesgue measure $\frac{dx}{2\pi}$, endowed with the standard scalar product $\Pa{f,g}:=\int_{S^1} f(x)\overline{g(x)}\ \frac{dx}{2\pi}$. An orthonormal basis is given by the set $\Pg{\fzRC}_{m\in\Z}$ with $\fzR(x):=x^{-m}$.

\bp
The continuous part of the dual space of $SU(1,1)$ is $$\hat{G}_C=\Pg{\RepC\ |\ j\in\Pg{0,\frac12},s=\frac12+iv, v\in\R^+}.$$

Given $\gsu\in SU(1,1)$, the unitary representation $\RepC(\gsu)$ is
\mmfunz{\RepC(\gsu)}{\domRC}{\domRC}
{\fzRbase(x)}{|\bar{\beta} x + \bar{\al}|^{-2s}
\Pt{\frac{\bar{\beta} x + \bar{\al}}{|\bar{\beta} x + \bar{\al}|}}^{2j} \fzRbase\Pt{\frac{\al x + \beta}{\bar{\beta} x + \bar{\al}}}}
with $\gsu^{-1}=\Mat{cc}{\al&\beta\\\bar{\beta}&\bar{\al}}$.

The Plancherel measure on $\hat{G}_C$ is $dP(j,\frac12+iv)=\begin{cases}
\frac{v}{2\pi}\Tanh(\pi v)\ dv\qquad j=0\\
\frac{v}{2\pi}\Ctanh(\pi v)\ dv\qquad j=\frac12.
\end{cases}$ where $dv$ is the Lebesgue measure on $\R$.
\ep
\brem
Notice that in this example the domain of the representation $\domRC$ does not depend on $j,s$.
\erem
\newcommand{\Ln}{\mathscr{L}_n}

Now we turn our attention to the description of principal discrete representations\footnote{There exist also the so-called complementary discrete representations, whose Plancherel measure is vanishing. Hence they do not contribute to the \GFT\ of a function defined on $SU(1,1)$. For details, see for instance \cite{sugiura}.}.

We first define the domain $\domRD$ of these representations $\RepD$: consider the space $\Ln$ of $L^2$ complex-valued functions on the unitary disc $D=\Pg{x\in\C\ |\ |x|< 1}$ with respect to the measure $d\mu^*(z)=\frac{2|n|-1}{\pi}(1-|z|^2)^{2n-2}\ dz$ where $dz$ is the Lebesgue measure. $\Ln$ is an Hilbert space if endowed with the scalar product $\Pa{f,g}:=\int_{D} f(z)\overline{g(z)}\ d\mu^*(z)$. Then define the space $\domRD$ with $n>0$ as the Hilbert space of holomorphic functions of $\Ln$, while $\domRD$ with $n<0$ is the Hilbert space of antiholomorphic functions of $\mathscr{L}_{-n}$. An orthonormal basis for $\domRD$ with $n>0$ is given by $\Pg{\fzRD}_{m\in\N}$ with $\fzRD(z)=\Pt{\frac{\Ga(2n+m)}{\Ga(2n)\Ga(m+1)}}^\frac12 z^m$ where $\Ga$ is the Gamma function. An orthonormal basis for $\domRD$ with $n<0$ is given by $\Pg{\fzRD}_{m\in\N}$ with $\fzRD(z)=\overline{\fzRbase_{-n}^m}(z)$.

\bp
The discrete part of the dual space of $SU(1,1)$ is $\hat G_D=\Pg{\RepD\ |\ n\in\frac12\Z, |n|\geq 1}$

Given $\gsu\in SU(1,1)$, the unitary representation $\RepD(\gsu)$ is
\mmfunz{\RepD(\gsu)}{\domRD}{\domRD}{\fzRbase (x)}{(\bar{\beta} x + \bar{\al})^{-2|n|}
\fzRbase \Pt{\frac{\al x + \beta}{\bar{\beta} x + \bar{\al}}}}
with $\gsu^{-1}=\Mat{cc}{\al&\beta\\\bar{\beta}&\bar{\al}}$.

The Plancherel measure on $\hat{G}_D$ is $dP(n)=\frac{2|n|-1}{4\pi}\cont(n)$, where $\cont$ is the counting measure.
\ep

\subsubsection{The kernel of the hypoelliptic heat equation}

In this section we compute the representation of differential operators $L_{X_i}$ with $i=1,2$ and give the explicit expression of the kernel of the hypoelliptic heat equation.

We first study the continuous representations $\dRepC{i}$, for both the families $j=0,\frac12$. Their actions on the basis $\Pg{\fzRC}_{m\in\Z}$ of $\domRC$ is
\bqn
\dRepC{1}\fzRC=\frac{s-m-j}{2}\fzRC[-1]+\frac{s+m+j}{2}\fzRC[+1],\nn
\dRepC{2}\fzRC=i\frac{s-m-j}{2}\fzRC[-1]-i\frac{s+m+j}{2}\fzRC[+1].
\eqnn
Hence $$\DHLC\fzRC=-\Pt{(m+j)^2+v^2+\frac{1}{4}}\fzRC.$$
Moreover, the set 
$\Pg{\fzRC}_{m\in\Z}$ is a complete set of eigenfunctions of norm one for the operator $\DHLC$.

\brem
Notice that the operators $\dRepC{i}$ are only defined on the space of $C^\infty$ vectors, i.e. the vectors $v\in\domRC$ such that the map $g\rightarrow [\RepC(g)]\ v$ is a $C^\infty$ mapping. This restriction is not crucial for the following treatment.
\erem

We now turn our attention to the discrete representations in both cases $n>0$ (holomorphic functions) and $n<0$ (antiholomorphic functions). Consider the discrete representation $\RepD$ of $SU(1,1)$ and let $\dRepD{i}$ be the representations of the differential operators $L_{X_i}$ with $i=1,2$. Their actions on the basis $\Pg{\fzRD}_{m\in\N}$ of $\domRD$ are
\bqn
\dRepD{1}\fzRD&=&\frac{\sqrt{(2|n|+m)(m+1)}}{2}\fzRD[+1]
-\frac{\sqrt{(2|n|+m-1)m}}{2}\fzRD[-1]\nn
\dRepD{2}\fzRD&=&-i\frac{\sqrt{(2|n|+m)(m+1)}}{2}\fzRD[+1]
-i\frac{\sqrt{(2|n|+m-1)m}}{2}\fzRD[-1]
\eqnn
Hence $\DHLD\fzRD=-\Pt{|n|+2m|n|+m^2}\fzRD$, thus the basis
$\Pg{\fzRD}_{m\in\N}$ is a complete set of eigenfunctions of norm one for the operator $\DHLD$.

We now compute the kernel of the hypoelliptic heat equation using Formula \r{e-fundsol-general}.
\bp
The kernel of the hypoelliptic heat equation on $(SL(2),\distr,\metr)$ is
\bqn
p_t(g)&=&\int_0^{+\infty} \frac{v}{2\pi}\Tanh(\pi v) \sum_{m\in\Z} e^{-t\Pt{m^2+v^2+\frac{1}{4}}}\Pa{\fzRC,\RepC[0](\gsu) \fzRC}\ \leb(v) + \nn &+&\int_0^{+\infty} \frac{v}{2\pi}\Ctanh(\pi v) \sum_{m\in\Z} e^{-t\Pt{m^2+m+v^2+\frac12}}\Pa{\fzRC,\RepC[\frac12](\gsu) \fzRC}\ \leb(v) + \nn
&+& \sum_{n\in\frac12\Z,\ |n|\geq 1} \frac{2|n|-1}{4\pi} \sum_{m\in\N}
e^{-t(|n|+2m|n|+m^2)} \Pa{\fzRD,\RepD(\gsu)\fzRD}.
\eqnl{eq-SL2-heat-exp}
where $\gsu=\Pi(g^{-1})\in SU(1,1)$.
\ep

\subsection{The hypoelliptic heat kernel on $\mot$}
\renewcommand{\Rep}{\Repbase^\lam}
\renewcommand{\domR}{\domRbase}
\renewcommand{\dRep}[1]{\dRepbase{#1}{\lam}}
\label{s-mot}

%%%%%%%%%%%%%%%%%%%%%%%%%%%%%%%%%%%%%%%%%%%%%%%%%%%%%%%%%%%%%%%%%%%%
Consider the group of rototranslations of the plane
\bqn
\mot=\Pg{\left(
\ba{ccc}
\cos(\al)&-\sin(\al)&x_1\\
\sin(\al)&\cos(\al)&x_2\\
0&0&1
\ea
\right)\ |\ \al\in\R/2\pi,\ x_i\in\R}
\eqnn
In the following we often denote an element of $\mot$ as $g=(\al,x_1,x_2)$.

A basis of the Lie algebra of $\mot$ is $\Pg{p_0,p_1,p_2}$, with
\bqn
p_0=\left(
\ba{ccc}
0&-1 &0 \\
1&0 &0 \\
0 &0 &0 
\ea
\right),\hspace{.5cm} 
p_1=\left(
\ba{ccc}
0&0 &1 \\
0&0 &0 \\
0 &0 &0 
\ea
\right),\hspace{.5cm} 
p_2=\left(
\ba{ccc}
0&0 &0 \\
0&0 &1 \\
0 &0 &0 
\ea
\right)
\eqn

We define a trivializable sub-Riemannian structure on $\mot$ as presented in Section \ref{ss-leftmanifold}: consider the two left-invariant vector fields $X_i(g)=g p_i$ with $i=0,1$ and define
\bqn
\distr(g)=\span{X_0(g),X_1(g)}~~~~~~\metr_g(X_i(g),X_j(g))=\de_{ij}.
\eqnn

The group $\mot$ is unimodular, hence the hypoelliptic Laplacian $\Dh$ has the following expression:
%%%%%%%%%%%%%%%%%%%%%%%%%%%%%%%%%%%%%%%%%%%%%%%%%%
\bqn
\llabel{eq-M2-heat}
\Dh\phi=\left(L_{X_0}^2+ L_{X_1}^2  \right) \phi
\eqn

\brem
As for the Heisenberg group, all left-invariant sub-Riemannian structures that one can define on $\mot$  are isometric.
\erem

In the following proposition we present the structure of the dual of $\mot$.
\bp
The dual space of $\mot$ is $\hat{G}=\Pg{\Rep\ |\ \lam\in\R^+}$.

Given $g=(\al,x_1,x_2)\in \mot$, the unitary representation $\Rep(g)$ is
\mmfunz{\Rep(g)}{\domR}{\domR}{\fzRbase(\th)}{e^{i\lam(x\cos(\th)-y\sin(\th))}\fzRbase(\th+\al),}
where the domain $\domR$ of the representation $\Rep(g)$ is $\domR=L^2(S^1,\C)$, the Hilbert space of $L^2$ functions on the circle $S^1\subset\R^2$ with respect to the Lebesgue measure $d\th$, endowed with the scalar product $<\fzRbase_1,\fzRbase_2>=\int_{S^1}\fzRbase_1(\th)\overline{\fzRbase_2(\th)}\ d\th$.

The Plancherel measure on $\hat G$ is $dP(\lam)=\lam d\lam$ where $d\lam$ is the Lebesgue measure on $\R$.
\ep
\brem
Notice that in this example the domain of the representation $\domR$ does not depend on $\lam$.
\erem

\subsubsection{The kernel of the hypoelliptic heat equation}

Consider the representations $\dRep{i}$ of the differential operators $L_{X_i}$  with $i=1,2$: they are\hyphenation{o-pe-ra-tors} operators on $\domR$, whose action on $\fzRbase\in\domR$ is (using formula \r{diff-rep})
\bqn
\Pq{\dRep{0} \fzRbase}(\th)=\frac{d\fzRbase(\th)}{d\th}\qquad
\Pq{\dRep{1} \fzRbase}(\th)=i\lam\cos(\th)\fzRbase(\th),\eqnn
hence
\bqn
\Pq{\DHL \fzRbase}(\th)=\frac{d^2\fzRbase(\th)}{d\th^2}- \lam^2\cos^2(\th)\fzRbase(\th).
\eqnn

We have to find a complete set of eigenfunctions of norm one for $\DHL$. An eigenfunction $\fzRbase$ with eigenvalue $E$ is a $2\pi$-periodic function satisfying the {\bf Mathieu's equation}
\bqn\llabel{eq-mathieu}
\frac{d^2\fzRbase}{dx^2}+(a-2q\cos(2x))\fzRbase=0
\eqn
with $a=-\frac{\lam^2}{2}-E$ and $q=\frac{\lam^2}{4}$. For details about Mathieu functions see for instance \cite[ch. 20]{mathieu}.

\brem
Notice that we consider only $2\pi$-periodic solutions of \r{eq-mathieu} since $\domR=L^2(S^1,\C)$.
\erem

There exists an ordered discrete set $\Pg{a_n(q)}_{n=0}^{+\infty}$ of distinct real numbers ($a_{n}<a_{n+1}$) such that the equation $\frac{d^2f}{dx^2}+(a_n-2q\cos(2x))f=0$ admits a unique even $2\pi$-periodic solution of norm one. This function $\ce_n(x,q)$ is called an {\bf even Mathieu function}.

Similarly, there exists an ordered discrete set $\Pg{b_n(q)}_{n=1}^{+\infty}$ of distinct real numbers ($b_{n}<b_{n+1}$) such that the equation $\frac{d^2f}{dx^2}+(b_n-2q\cos(2x))f=0$ admits a unique odd  $2\pi$-periodic solution of norm 1. This function $\se_n(x,q)$ is called an {\bf odd Mathieu function}.

The set $\mathcal{B}^\lam:= \Pg{\ce_n\Pt{x,\frac{\lam^2}{4}}}_{n=0}^{+\infty}\cup\Pg{\se_n\Pt{x,\frac{\lam^2}{4}}}_{n=1}^{+\infty}$ is a complete set of $2\pi$-periodic eigenfunctions of norm one for the operator $\DHL$. The eigenvalue for $\ce_n\Pt{x,\frac{\lam^2}{4}}$ is $a^\lam_n:=-\frac{\lam^2}{2}-a_n\Pt{\frac{\lam^2}{4}}.$ The eigenvalue for $\se_n\Pt{x,\frac{\lam^2}{4}}$ is $b^\lam_n:=-\frac{\lam^2}{2}-b_n\Pt{\frac{\lam^2}{4}}.$ 

We can now compute the explicit expression of the hypoelliptic kernel on $\mot$.
\bp
The kernel of the hypoelliptic heat equation on $(\mot,\distr,\metr)$ is
\bqn
\llabel{eq-M2-heat-exp}
p_t(g)= \int_0^{+\infty}\lam\ d\lam\Pt{
\sum_{n=0}^{+\infty} e^{a_n^\lam t} <\ce_n(\th),\Rep(g)\ce_n(\th)>+
\sum_{n=1}^{+\infty} e^{b_n^\lam t} <\se_n(\th),\Rep(g)\se(\th)>}
\eqn
The function \r{eq-M2-heat-exp} is real for all $t>0$.
\ep
\proof The formula \r{eq-M2-heat-exp} is given by writing the formula \r{e-fundsol-general} in the $\mot$ case.

We have to prove that $p_t(g)$ is real: we claim that $<\ce_n,\Rep(g)\ce_n>$ is real. In fact, write the scalar product with $g=(\al,x,y)$:
$$<\ce_n,\Rep(g)\ce_n>=\int_0^{2\pi}e^{i\lam(x\cos(\th)-y\sin(\th))}\ce_n(\th)\ce_n(\th+\al).$$ Its imaginary part is $\int_0^{2\pi}\sin\Pt{\lam(x\cos(\th)-y\sin(\th))}\ce_n(\th)\ce_n(\th+\al)$. Its integrand function assumes opposite values in $\th$ and $\th+\pi$: indeed $$\sin\Pt{\lam(x\cos(\th+\pi)-y\sin(\th+\pi))}=
\sin\Pt{\lam(-x\cos(\th)+y\sin(\th))}=-\sin\Pt{\lam(+x\cos(\th)-y\sin(\th))},$$ while $\ce_n(\th+\pi)=(-1)^n\ce_n(\th)$ as a property of Mathieu functions. Thus, the integral over $[0,2\pi]$ is null. With similar observations it is possible to prove that $<\se_n(\th),\Rep(g)\se_n(\th)>$ is real.

Thus, $p_t(g)$ is an integral of a sum of real functions, hence it is real.
\qed\\\\

{\bf Acknowledgments.} The authors are grateful to Giovanna Citti and to Fulvio Ricci for helpful discussions.

\end{document}